\documentclass[12pt,noindent]{article}

\usepackage{amsmath}
\usepackage{amsthm}
\usepackage{amssymb}
\usepackage{amsfonts}
\usepackage{amsxtra}                     
\usepackage{graphicx}
\usepackage{latexsym}
\usepackage{verbatim}
\include{header}

\theoremstyle{plain}
\newtheorem{theorem}{Theorem}

\newtheorem{lemma}[theorem]{Lemma}

\newtheorem{definition}[theorem]{Definition}

\newtheorem{claim}[theorem]{Claim}

\newtheorem{remark}[theorem]{Remark}

\newcommand{\thmref}[1]{Theorem~\ref{thm:#1}} 
\newcommand{\lemref}[1]{Lemma~\ref{lem:#1}} 
\newcommand{\remref}[1]{Remark~\ref{rem:#1}} 
\newcommand{\claimref}[1]{Claim~\ref{claim:#1}} 
\newcommand{\secref}[1]{Section~\ref{sec:#1}} 
\newcommand{\subref}[1]{Subsection~\ref{sub:#1}} 
\newcommand{\eqnref}[1]{(\ref{eq:#1})} 
\newcommand{\figref}[1]{Figure~\ref{fig:#1}} 

\newcommand\ignore[1]{}

\newcommand\mypar[1]{\noindent\par{\sc #1}~}

\def\R{\mathbb{R}} 
\def\N{\mathbb{N}} 

\def\periodeq{\mbox{.}}
\def\commaeq{\mbox{,\,}}

\def\suchthat{\,\mbox{:}\,}

\renewcommand{\Pr}[1]{\mathbb{P}\left(#1\right)} 
\newcommand{\Ex}[1]{\mathbb{E}\left[#1\right]} 
\def\eqdist{=^d}   

\newcommand{\bigoh}[1]{O\left(#1\right)}

\newcommand{\ohmega}[1]{\Omega\left(#1\right)}
\newcommand{\theita}[1]{\Theta\left(#1\right)}

\def\sE{\mathcal{E}}

\def\sR{\mathcal{R}}

\def\Ve{{V}} 
\def\Ed{{E}} 

\newcommand\QED{\ifhmode\allowbreak\else\nobreak\fi
\quad\nobreak$\Box$\medbreak}
\newcommand{\proofstart}{\par\noindent\sl Proof:\rm\enspace}
\newcommand{\proofend}{\QED\par}
\renewenvironment{proof}{\proofstart}{\proofend}

\def\eps{\epsilon}

\def\periodeq{\mbox{.}}
\def\commaeq{\mbox{,}}

\def\E{\sE}

\def\Efin{\sE^{\text{fin}}}

\def\deg{\indeg}
\newcommand{\smL}{\mbox{\sf sm}_L}
\newcommand{\lgL}{\mbox{\sf lg}_L}
\def\hati{I_i}
\newcommand{\estr}{\epsilon}
\newcommand{\tree}{\text{T}}
\newcommand{\ftree}{\mbox{T}_{\infty}}

\newcommand{\glue}{\mbox{\sc Glue}(\tree,v,k)}
\def\KR{\mbox{\sf GN}}

\def\proper{parent-closed}

\newcommand{\W}{\mathcal{W} }
\newcommand{\T}{\mathcal{T}}
\renewcommand{\deg}{\mbox{d}}
\newcommand{\s}{\mathcal{S}}
\newcommand{\p}{\mathcal{P}}
\newcommand{\B}{\mathcal{B}}

\def\Nstr{\mathbb{N}^*}

\def\nat{\N}
\def\natstr{\Nstr}
\def\real{\mathbb{R}}
\newcommand{\clref}[1]{Claim~\ref{cl:#1}}
\newcommand{\pr}[1]{\Pr{#1}}

\newtheorem{fact}{Fact}

\numberwithin{theorem}{section}        
\numberwithin{equation}{section} 

\begin{document}
\title{Connectivity transitions in networks with super-linear preferential attachment}

\author{Roberto Oliveira\thanks{IBM T. J. Watson Research Laboratory, Yorktown Heights, NY 10598. E-mail: riolivei@us.ibm.com. Work done while the author was a Ph.D. student at the Courant Institute, New York University, supported by a fellowship from
CNPq, Brazil.} \and Joel Spencer \thanks{Courant Institute, New
York University, New York, NY 10012. E-mail:
spencer@cims.nyu.edu}} \maketitle \abstract{We analyze an evolving
network model of Krapivsky and Redner in which new nodes arrive
sequentially, each connecting to a previously existing node $b$
with probability proportional to the $p$-th power of the in-degree
of $b$.  We restrict to the super-linear case $p>1$.  When
$1+\frac{1}{k} < p < 1 + \frac{1}{k-1}$ the structure of the final
countable tree is determined.  There is a finite tree $\mbox{T}$
with distinguished $v$ (which has a limiting distribution) on
which is ``glued" a specific infinite tree.  $v$ has an infinite
number of children, an infinite number of which have $k-1$
children, and there are only a finite number of nodes (possibly
only $v$) with $k$ or more children.  Our basic technique is to
embed the discrete process in a continuous time process using
exponential random variables, a technique that has previously been
employed in the study of balls-in-bins processes with feedback. }

\bibliographystyle{plain}       


\section{Introduction}\label{sec:intro}
In some important examples of growing networks, such as the World
Wide Web or the scientific citation network, one can interpret the
fact that a given node has high in-degree as indicative that node
is ``popular''. For instance, popular papers are the ones more
often cited more by other works, and popular Web pages receive
more links than less popular ones. A consequence of differences in
popularity is that a node with high in-degree has more propensity
to receive further edges as the network evolves than an unpopular
node with low in-degree. In other words, the more popular a node
is, the more visible it is to the community that creates the
network and/or interacts through it, and high
visibility makes future increases in popularity more likely.\\

Barab\' asi and Albert \cite{Barabasi99} incorporated this
so-called {\em preferential attachment} phenomenon into a
generative model for these and other networks. In this model,
nodes arrive at the network one at a time, and direct a fixed
number $m$ of edges to previously existing nodes that are chosen
with probabilities proportional to their in-degrees. It is quite
remarkable that this simple model already replicates many
non-trivial features of the above networks, such as power-law
degree
distributions, small diameter and high resistance to random failures, as argued non-rigorously by physicists (see \cite{Albert99,AlbertSurvey} and references therein) and later proven rigorously by mathematicians \cite{Bollobas00,Bollobas03,Bollobas04,BollobasSurvey}.\\

The success of the Barab\'asi-Albert model has also inspired many
different variants. The models in \cite{Cooper03,Bollobas03_2}
permit that the power-law exponent of the degree distribution be
adjusted to fit real-world data. Other models \cite{Bianconi01}
feature preferential attachment that is dictated both by node
fitness and popularity. This work is dedicated to yet another kind
of variant of the model of \cite{Barabasi99}, one in which the
{\em strength} of preferential attachment
can be varied.\\

This model was proposed and studied by Krapivsky and Redner
\cite{Krapivsky02} and independently by Drinea, Mitzenmacher and
Enachescu \cite{DrineaEM}. It differs from the Barab\'asi-Albert
network in that each incoming node chooses a pre-existing vertex
to link to with probability proportional to a fixed function $f$
(the {\em attachment kernel}) of the degree of that
vertex\footnote{As in the original Barab\' asi-Albert model, One
could also consider a similar model in which each incoming node
creates a fixed number $m$ of new edges, but we will only consider
the case $m=1$ in this paper}. While the Barab\'asi-Albert model
is recovered by setting $f(x)=x$, we will be mostly concerned with
kernels of the form $f(x)\sim x^p$ with $p>1$ thought of as a
tunable parameter; this is referred to in \cite{Krapivsky02} as
the {\em super-linear case}. One of the many remarkable
non-rigorous results about this so-called $\KR$ (Growing Network)
model is that it undergoes an infinite sequence of {\em
connectivity transitions} at $p=p_k\equiv 1 + 1/k$,
$k=1,2,3,\dots$. By this it is meant that for $p>p_{k}$ the $\KR$
process has only finitely many vertices that receive more than $k$
links, whereas for $p\leq p_k$ the number of such vertices is
infinite. Another way of stating this property is the following:
the smallest integer $k$ for which $p>p_k=1+1/k$ is also the
smallest number $k$ for which only finitely
many nodes ever reach in-degree $k$.\\

The connectivity transitions are both mathematically intriguing
and physically interesting. The fact that the $p>1$, $p=1$ and
(conjecturally) $p<1$ cases of the model are very different leads
the authors of \cite{MendesSurvey} to suggest that so-called {\em
self-organized criticality} is at work in networks with power-law
degree distributions (the $p=1$ case). It was also noted elsewhere
\cite{AlbertSurvey} that the condensation regime of the fitness
model of Bianconi and Barab\'asi \cite{Bianconi01} has
qualitatively similar behavior to the super-linear $\KR$;
\cite{AlbertSurvey} even suggests that a direct connection between
the two models could exist. There is also some modelling interest
in the connectivity transitions, since networks in which
preferential attachment is very strong (conceivably even some
parts of the World-Wide Web) should exhibit behavior that is qualitatively similar to the $\KR$ model in the super-linear regime.\\

Despite the striking characteristics, we do not know of any
rigorous work on the $\KR$  model to the present date. A modified
model was addressed in independent work by Chung, Handjani and
Jungreis \cite{Chung03}. In their process, an attachment kernel is
still present, but at each time step either a new vertex and a new
edge are added with probability $0<q<1$, or only a new edge is
added with probability $1-q$. This modified model does exhibit
connectivity transitions in the sense of \cite{Krapivsky02}, but
it is not clear how to deduce the analogous results for the
original $\KR$ model from
the techniques in \cite{Chung03}.\\

In this paper we attempt to give a rigorous description of the
super-linear $\KR$ process in the large-time limit. Our rigorous
results imply the existence of connectivity transitions, but they
also go beyond that. The first result we prove is the following.

\begin{theorem}\label{thm:KRintro_main1} Let $\{T_m\}_{m\geq 1}$ be the $\KR$ process with attachment kernel $f(x)= (x+1)^p$ (defined in \secref{KRintro_defin}). Also let $\ftree$ be the increasing limit of the $\{T_m\}_{m\geq 1}$ process, and assume that $p>p_k = 1 + 1/k$. Then with probability $1$ all but finitely many nodes of $\ftree$ have less than $k$ descendants (cf. the definition in \subref{treeterm}).\end{theorem}

A vertex of $\ftree$ with in-degree larger than or equal to $k$
necessarily has at least $k$ descendants. For this reason,
\thmref{KRintro_main1} implies that for $p>p_k$, only finitely
many vertices of $\ftree$ have in-degree $\geq k$. As a result,
the number of vertices in $\tree_m$ with in-degree bigger than $k$
is bounded as $m\to +\infty$. This shows that
\thmref{KRintro_main1} implies the non-rigorous ``$p>p_k$" result
of \cite{Krapivsky02}, and is in fact stronger than it. Similarly,
\thmref{KRintro_main2} below implies the $p\leq p_k$ case of
Krapivsky and Redner's result.

\begin{theorem}\label{thm:KRintro_main2}Let $f$, $p>1$ and $\ftree$ be as in \thmref{KRintro_main1}, and let $k=k_p$ be the smallest positive integer for which
$p>p_k= 1+1/k$. Consider the construction $\glue$ defined in
\subref{treeterm}. Then the set of values (up to isomorphism) that
$\ftree$ attains with positive probability is precisely the set of
all trees that can be obtained by choosing a finite rooted tree
$T$, a distinguished vertex $v\in T$ and setting
$\ftree=\glue$.\end{theorem}

\thmref{KRintro_main2} completely describes (up to isomorphisms)
the limit set of the $\KR$ process in the large-time limit. In
particular, it also implies that if $p\leq p_k$, the number of
vertices of in-degree $\geq k$ in $\tree_m$ diverges as $m\to
+\infty$. This differs from the original claim in
\cite{Krapivsky02}, in which the authors argue that the
expectation of the number of vertices of degree $\ell \geq k$
diverges at a certain rate. While we have nothing to say about
this rate, \thmref{KRintro_main2} is stronger than the claim of
\cite{Krapivsky02} in that divergence of the expected number is
implied, but does not imply, almost sure divergence. Moreover, our
description of the structure of $\ftree$ is new. Finally, we note
that there is nothing special about the choice of $f(x)=(x+1)^p$
as our superlinear kernel. In fact, the proof of both theorems
will make it clear that it suffices to assume that $f(x)>0$ for
all $x$ and that $f(x)=\theita{x^p}$ for $x\gg 1$, with only minor modifications in our arguments.\\

We now briefly outline our proof techniques. On a high level, we
rely strongly on the similarity pointed out by Drinea, Frieze and
Mitzenmacher \cite{Drinea02} between the $\KR$ process and {\em
balls-in-bins models with feedback}. The latter model describes
the evolution of a system with a fixed number of bins at which
balls are thrown. A ball arrives at each discrete time step and
chooses a bin to go into with probability proportional to a fixed
function $f$ (that we call the feedback function) of the number of
balls currently in that bin. This model can also be viewed as a
static variant of the $\KR$ process in which new edges are
repeatedly added but without the creation of any new nodes/bins .
This analogy permits that a certain technique applied to the study
of balls-in-bins problems
\cite{Khanin01,MitzenmacherOS04,Spencer??} is adapted to the $\KR$
process. It consists of building a continuous-time process out of
exponential random variables and showing that it embeds the
original discrete-time process. For this reason we call this
construction the {\em exponential embedding}. The $\KR$ version of
the exponential embedding is essential to the construction and
analysis of the infinite tree limit $\ftree$,
and we view it as an important part of our paper's contribution.\\

The remainder of the paper is organized as follows. In
\secref{KRintro_prelim} we introduce our notation and review a few
basic concepts. We formally define the $\KR$ process in
\secref{KRintro_defin}, starting with its original definition in
\cite{Krapivsky02,DrineaEM}, and then describing a useful labelled
version of it. \secref{KRintro_expemb} introduces the exponential
embedding technique. We begin with a review of the simpler
balls-in-bins case, then move on to the construction of the
embedding of the $\KR$ process for general attachment kernels. We
then employ the embedding to show that so-called ``explosive
kernels'' give rise to $\KR$ processes for which $\tree_m\to
\ftree$ in finite time under the exponential embedding. This
section ends with some lemmas on sums of exponential random
variables that will be useful later on. \thmref{KRintro_main1} is
proven in the subsequent \secref{KRmain1}. The section starts with
weaker results that intuitively pave the way for the actual proof
of the Theorem, which relies on a careful consideration of the
time of the birth of the $k$th descendant of a given node in the
exponential embedding setting. In \secref{KRmain2} we prove
\thmref{KRintro_main2}, relying on \thmref{KRintro_main1} and on
the techniques developed in the previous sections. We discuss some
consequences of our main
theorems and some related open questions in the Conclusion (\secref{KRextra}). The Appendix contains the proofs of some technical results.\\

\mypar{Acknowledgements.} We thank Eleni Drinea and Michael
Mitzenmacher for bringing this problem to our attention and for
useful discussions. We also thank the anonymous referees for
pointing out several typos and making suggestions that greatly
improved our presentation.

\section{Preliminaries}\label{sec:KRintro_prelim}

\subsection{Probabilistic ingredients}\label{sub:probing}

We briefly remind the reader of some basic probabilistic concepts
and tools, while also fixing some notation.\\

\mypar{Distributions.} We say that two random variables $X$, $Y$
taking values on the same set $U$ {\em have the same distribution}
(or {\em are identical in law}) if for all measurable subsets
$A\subseteq U$ $\Pr{X\in A} = \Pr{Y\in A}$. This will be
symbolically represented by $X\eqdist Y$.\\

\mypar{The exponential distribution.} A random variable $X$ is
said to be {\em exponentially-distributed with rate $\lambda>0$}
if $X$ almost surely takes values on the positive reals and
$$\Pr{X>t} = e^{-\lambda t}\;\; (t\geq 0)$$
We denote this property by $X\eqdist \exp(\lambda)$. The shorthand
$\exp(\lambda)$ will also denote a generic
exponentially-distributed random variable with rate $\lambda$. We
list below some elementary but extremely useful properties of
those random variables.
\begin{enumerate}\item {\em Lack of memory.} Let $X\eqdist \exp(\lambda)$
and $Z\geq 0$ be independent from $X$. The distribution of $X-Z$
conditioned on $X>Z$ is still equal to $\exp(\lambda)$. \item {\em
Minimum property.} Let $\{X_i\eqdist \exp(\lambda_i)\}_{i=1}^{m}$
be independent. Then $X_{\min}\equiv \min_{1\leq i\leq
m}X_i\eqdist \exp(\lambda_1+\lambda_2+\dots \lambda_m)$. Moreover,
for all $1\leq i\leq m$,
$$\Pr{X_i = X_{\min}} = \frac{\lambda_i}{\lambda_1+\lambda_2+\dots
\lambda_m}\mbox{ .}$$  \item {\em Multiplication property.} If
$X\eqdist\exp(\lambda)$ and $\eta>0$ is a fixed number, $\eta X
\eqdist\exp(\lambda/\eta)$.\end{enumerate}

 \mypar{The Borel-Cantelli Lemma.} Let $\{A_n\}_{n\in
N}$ be a sequence of events in some fixed probability space, with
$N$ a countable set. The event ``$A_n$ {\em infinitely often}
{$(n\in N)$}" (or ``$A_n$ {\em i.o.} {$(n\in N)$}") contains all
outcomes that belong to an infinite number of the events $A_n$.
The {\em Borel-Cantelli Lemma} states that
$$\sum_{n\in N}\Pr{A_n} <+\infty \Rightarrow \Pr{A_n\text{ i.o. }(n\in N)}=0$$
and
$$\sum_{n\in N}\Pr{A_n} =+\infty \text{  and  }\{A_n\}_{n\in N}\text{ independent }\Rightarrow \Pr{A_n\text{ i.o. }(n\in N)}=1\mbox{ .}$$

\mypar{Discrete-time Markov Chains.} A (discrete-time) Markov
chain on the countable set $\Omega$ is specified by {\em
transition probabilities} $\Pi:\Omega\times\Omega \to [0,1]$ and a
{\em initial condition} $X_0\in \Omega$ (possibly
non-deterministic). The recipe

$$\Pr{\forall 0\leq i\leq t \;\; X_i=\omega_i} = \Pr{X_0=\omega_0}\prod_{i=1}^{t}\Pi(\omega_{i-1},\omega_i)$$
defines the distribution of a sequence $\{X_i\}_{i=0}^{+\infty}$
of $\Omega$-valued random variables.

\subsection{Tree terminology}\label{sub:treeterm}

\mypar{Trees.} All trees are rooted and have their edges directed
towards the root. No loops or parallel edges are allowed.
Given vertices $a,b$ in a tree $T$, the existence of the oriented edge $(a,b)$ will be indicated by saying that $a$ {\em is a child of }$b$, or that $b$ {\em is $a$'s parent}, or that $a$ {\em links to }$b$. With this terminology, the (in-)degree $\deg_T(b)$ of $b$ in $T$ is the number of its children. If $r$ is a node of $T$, the subtree $T_r$ of $T$ rooted at $r$ is the tree with root $r$, together with $r$'s children, the children of those children, and so on. The nodes in $T_r\backslash \{r\}$ are referred to as the descendants of $r$, and $r$ is said to be {\em $k$-fertile} in $T$ if it has $k$ or more descendants. \thmref{KRintro_main1} consists of showing that for $p>p_k$, only finitely many nodes in $\ftree$ are $k$-fertile.\\

\mypar{The $\mbox{\sc 'Glue'}$~construction.} Given a finite
(rooted, oriented) tree $\tree$, a distinguished node $v$ of
$\tree$ and an integer $k\geq 1$, we define $\glue$ as follows.
For each finite (rooted, oriented) tree $S$ on $k$ or less nodes,
take countably many copies $\{S_i\}_{i\geq 1}$ of $S$. $\glue$ is
the union of $\tree$
with all the trees $S_i$ as above, with the addition of edges from the root of each one of the $S_i$'s to $v$.\\

As a simple example suppose $\tree$ consists of a single node (the
root $v$) and $k=2$.  Then in $\glue$ the root has a countably
infinite number of children.  Infinitely many of these children
are childless and infinitely many of these children have precisely
one child and none of them have more than one child.  Further, all
grandchildren
of the root are childless.\\

\begin{figure}[t]
\centering \includegraphics{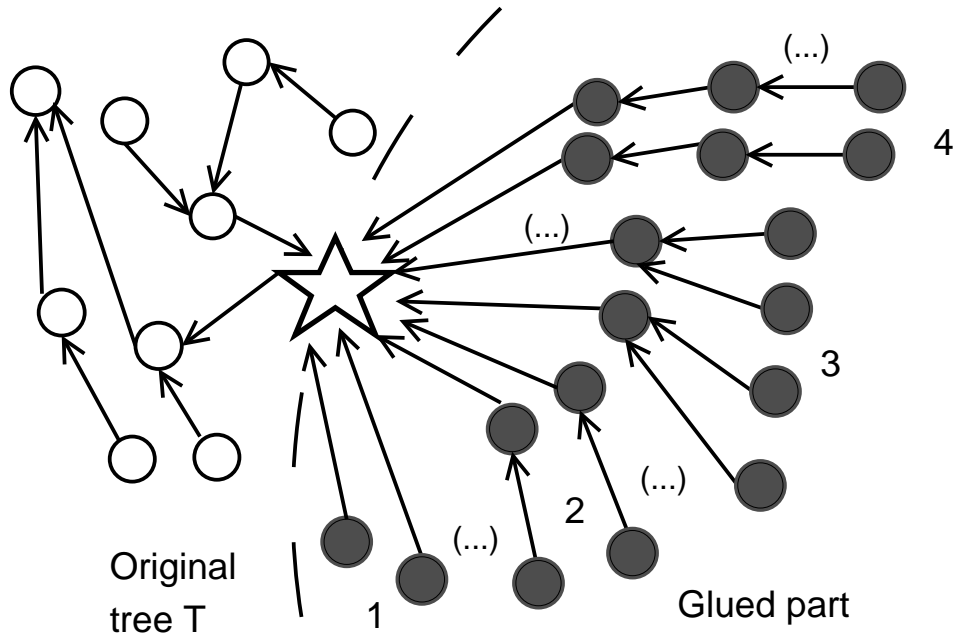}\caption{\normalsize An
example of $\glue$ for $k=3$.}\label{fig:glue}
\end{figure}

A more complex example of $\glue$, now with $k=3$, is portrayed in
\figref{glue}. The starred node is $v$, and the finite tree
$\tree$ lies to the left of the dashed line. The countably many
copies of the four rooted trees on $3$ or less
vertices (numbered $1$, $2$, $3$ and $4$ in the Figure) appear to the right of the line, and are all connected to $v$ by their roots.\\

Our task in proving \thmref{KRintro_main2} will be to show that
with probability $1$ there exist $\tree$ and $v$ as above with
$\ftree = \glue$, and that all such $\glue$ occur as values of
$\ftree$ with some positive probability.

\subsection{Labelled trees and \proper~sets}\label{sub:labelsmil}

\mypar{Labels.} It will be convenient for us to label the vertices
of trees. For our purposes, a label is a (possibly empty) sequence
of elements of the set $\nat=\{1,2,3,\dots\}$ of positive
integers. The empty sequence is denoted by $\estr$, and all other
sequences $a=\{a_i\}_{i=1}^m\subset \nat$ (with $m\geq 1$) will be
represented by
$a=a_1a_2\dots a_m$. Moreover, we call $a_1\dots a_{m-1}$ (the sequence $a$ without its last element) the {\em parent sequence} of $a$. The set of all labels will be denoted by $\natstr$.\\

\mypar{Labelling trees.} A labelling of a finite tree $T$ is an
assignment of labels to the vertices of $T$ that obeys two rules.

\begin{itemize} \item the label of the root of $T$ is the empty sequence $\estr$;
\item if vertex $v$ has degree $d$ and is labelled by the sequence
$v_1\dots v_m$, its children will receive labels $v_1\dots v_m i$,
$1\leq i\leq d$.\end{itemize}

The second rule implies that the label of a vertex $v$'s parent in
$T$ is the parent sequence of the label of $v$.\\

\mypar{Parent-closed subsets.} A subset $A\subset \natstr$ is said
to be {\em \proper~}if it is non-empty and for all non-empty
sequences $a\in A$ the parent sequence of $a$ is also in $A$. Any
\proper~$A$ corresponds to a finite tree with vertex set $A$ and
edges from each $a\in A\backslash\{\estr\}$ to $a$'s parent.
Conversely, given a tree $T$, the labelling procedure above
provides a proper set $A=A(T)$ that corresponds to tree $T$. This
set $A(T)$ is not uniquely defined, but this will not keep us from
representing finite trees by finite \proper~$A\subset \natstr$ in
what follows. For this reason, we will often apply tree
terminology to \proper~ $A\subset\natstr$, speaking for instance
of the degree $\deg_A(a)$ of an element $a\in A$. We also observe
that the potential descendants of $a=a_1\dots a_m\in \natstr$ are
obtained by adjoining the terms of another sequence $b=b_1\dots
b_n$ to $a$, thus forming the {\em concatenation} $ab\equiv
a_1\dots a_m b_1\dots b_n$. Finally, we define for convenience

$$\Efin \equiv \{A\subset \natstr\; : \; A\text{ \proper~ and
finite}\} \periodeq $$

\section{Definition of the $\KR$ process}\label{sec:KRintro_defin}
\subsection{The standard definition}\label{sub:KRintro_def}

The $\KR[f]$ process is defined in terms of an {\em attachment
kernel}, that is, a function $f:\N\cup\{0\}\to\R^+$. The process
evolves in discrete time $m=0,1,2,3,\dots$; its state at time
$m\geq 0$ is a (rooted, oriented) tree
$\tree_{m}=(\Ve_{m},\Ed_{m})$ with vertex set $\Ve_{m}$ and edge
set $\Ed_{m}$. Initially, $\tree_{0}$ contains a single {\em root}
node and no edges. At each time $m>0$ the tree is updated by the
addition of a new node ($\Ve_{m} = \Ve_{m-1} \cup \{v_m\}$) and a
new edge ($\Ed_{m}=\Ed_{m-1}\cup \{v_mw_m\}$), where $w_m$ is
chosen according to the following probability distribution:

$$\forall w\in
\Ve_{m-1}\; \Pr{w_m = w\mid \tree_{m-1}} =
\frac{f(\deg_{\tree_{m-1}}(w))}{\sum_{v\in \Ve_{m-1}}
f(\deg_{\tree_{m-1}}(v))}\mbox{ .}$$ These definitions already
specify the process completely as a finite-tree-valued Markov
Chain.

\subsection{The labelled process}\label{sub:KRintro_label}

In the labelled $\KR[f]$ process, we start by labelling the root
(and unique element) of $\tree_0$ by the empty sequence $\estr$.
At subsequent times $m\geq 1$, assume that the incoming node $v_m$
links to a node $w_m$ that is labelled by the sequence $a_1\dots
a_n$, and that $v_m$ is the $\ell$th node to link to $w_m$. Then
the label of $v_m$ is defined to be $a_1\dots a_n \ell$, i.e. the
sequence corresponding to $v_m$'s parent $w_m$, with a new number
$\ell$ added to it.\\

This recursive labelling obeys the definition of a labelling of a
tree given in \subref{labelsmil}, and provides an alternative
description of the process as a Markov Chain on $\Efin$, as
defined in \subref{labelsmil}. The transition probabilities of the
$\KR$ process on $\Efin$ are:
\begin{eqnarray}\label{eq:transition}\Pi(A,B) & = & \frac{f(\deg_A(a))}{\sum_{b\in A}f(\deg_A(b))}\commaeq \mbox{ if }\exists a\in A\suchthat B = A\cup \{a(\deg_A(a)+1)\}\\ \nonumber
& = & 0 \text{   otherwise}\end{eqnarray} and its initial state is
$\tree_0 = \{\estr\}$. We note  in passing that the limit $\ftree
= \bigcup_{m\geq 0}\tree_m$ of the $\KR$ process takes values in
the uncountable set
$$\E\equiv \{A\subset\Nstr: A\text{~\proper~and non-empty}\}\periodeq$$
$\E$ a closed subset of the topological space $2^{\Nstr}$ (with
the product topology). We will refrain from explicitly considering
measurability questions related to $\ftree$ and $\E$ in what
follows, since all such problems can be addressed in a rather
straightforward manner.

\section{Exponential embedding}\label{sec:KRintro_expemb}
Our aim in the present Section is to present the special
construction of the labelled $\KR$ process that we alluded to in
the Introduction. We will show how one can explicitly embed the
process in continuous time by employing sequences of independent
exponential
random variables. Although perhaps complicated at first sight, this embedding will prove to be fundamental to our analysis, with the independence of the involved random variables playing a key role in most of our computations.\\

\subsection{The balls-in-bins case}\label{sub:BinB}

Davis \cite{Davis90} applied the elementary properties of
exponential random variables to the study of Reinforced Random
Walks in a very interesting way. His method was later adapted by
Khanin and Khanin \cite{Khanin01} to the balls-in-bins setting. We
present this latter use of exponential random variables below
(which was also rediscovered by  Spencer and Wormald
\cite{Spencer??}) as a preparation
for the more difficult $\KR$ case.\\

Consider independent random variables $\{X_j,Y_j\eqdist
\exp(f(j))\}_{j\in\nat\cup\{0\}}$ and define, for $t\geq 0$
\begin{eqnarray*}N(t) & \equiv & \sup\{n\in\nat\cup\{0\} \mid \sum_{i=0}^{n-1}X_i\leq
t\}\commaeq\\
M(t) & \equiv & \sup\{m\in\nat\cup\{0\} \mid
\sum_{j=0}^{m-1}Y_j\leq t\}\periodeq\end{eqnarray*} We interpret
the times $\sum_{i=0}^{n-1} X_i$ and $\sum_{j=0}^{m-1} Y_j$ as the
times when $N(\cdot)$ and $M(\cdot)$ receive their $n$-th and
$m$-th ``hits", respectively. We now fix some $t\geq 0$ and $n,m
\in \nat\cup\{0\}$ and define the event
$$A^t_{n,m}\equiv \{N(t)=n,M(t)=m\} = \left\{\sum_{i=0}^{n-1} X_{i}\leq t <\sum_{i=0}^{n} X_{i}, \sum_{j=0}^{m-1} Y_{j}\leq t
<\sum_{j=0}^{m}Y_{j}\right\}\periodeq$$ What is the probability
that the $N(\cdot)$ process is the first one to receive a hit
after time $t$, conditioned on $A^t_{n,m}$? This probability can
be written as
$$\Pr{\left.\sum_{i=0}^{n} X_{i} < \sum_{j=0}^{m} Y_{j}\,\right|\,A^t_{n,m} }\periodeq$$
 If we further condition on $\sum_{i=0}^{n-1}X_i = s_1\leq
t$ and $\sum_{j=0}^{m-1}Y_j = s_2\leq t$, we can write this
probability as
$$\Pr{X_{n}+s_1 < Y_{m} + s_2\mid X_{n}\geq t-s_1, Y_m\geq t-s_2}$$
The lack-of-memory property of exponentials implies that under the
conditioning event above $X_n-t+s_1 \eqdist \exp(f(n))$ and
$Y_m-t+s_2\eqdist \exp(f(m))$. The minimum property then implies
\begin{multline*}\Pr{X_{n}+s_1 < Y_{m} + s_2\mid X_{n}\geq t-s_1, Y_m\geq t-s_2} =  \\ = \Pr{\exp(f(n))+t<\exp(f(m))+t} = \frac{f(n)}{f(n)+f(m)}\periodeq\end{multline*}
Since this holds for all $0\leq s_1,s_2\leq t$, we have in fact
proven that
$$\Pr{\left.\sum_{i=0}^{n} X_{i} < \sum_{j=0}^{m} Y_{j}\,\right|\,A^t_{n,m} } =  \frac{f(n)}{f(n)+f(m)}\periodeq$$
We thus arrive at a surprising conclusion.

\begin{fact}[Exponential embedding for balls-in-bins, \cite{Davis90,Spencer??}] Consider the balls-in-bins process \cite{Drinea02} with two bins and
feedback function $f$, i.e. the discrete Markov Chain that evolves
from state $(n,m)\in (\nat\cup\{0\})^2$ to state $(n+1,m)$ with
probability $\frac{f(n)}{f(n)+f(m)}$ and from $(n,m)$ to $(n,m+1)$
with probability $\frac{f(m)}{f(n)+f(m)}$. It then holds that the
joint hit counts of the $(N(\cdot),M(\cdot))$ processes up to the
(possibly finite) time when either one becomes infinite is
identical in law to the balls-in-bins process with feedback
function $f$ started from $(0,0)$. That is, the balls-in-bins
process is embedded in the continuous time $(N(\cdot),M(\cdot))$
process, with $X_j$ (respectively $Y_j$) parameterizing the time
between the arrivals of the $j$-th and $(j+1)$-th balls at the
first (resp. second) bin.\end{fact} Many non-trivial results that
do not have direct combinatorial proofs can be deduced from the
above construction. This method seems to be especially powerful in
the case when either $N$ or $M$ reaches an infinite value in
finite time. The reader is directed to
\cite{Khanin01,MitzenmacherOS04,Spencer??} for many examples of
applications of the exponential embedding. We will now show how we
can adapt this technique to our present context.

\subsection{Exponential embedding of the $\KR$ process}\label{sub:expemb} As pointed out in the introduction, a
balls-in-bins process with feedback function $f$ is very similar
to a $\KR$ process with attachment kernel $f$ to which only new
edges (and no new vertices) are added. Conversely, one may think
of a $\KR$ process as a balls-in-bins process in which each new
ball also creates a corresponding bin. This analogy was exploited
in \cite{Chung03}, in which a variant of the original $\KR$
process was modelled as an ``infinite P\'olya Urn process" for the
purposes of studying the degree sequence. We take this analogy
further by adapting the exponential embedding technique to
the labelled $\KR$ process as defined in \subref{KRintro_label}.\\

Our construction starts from an independent sequence
$\{X(a,j)\eqdist\exp(f(j))\mid a\in\natstr,j\in\nat\cup\{0\}\}$ of
random variables. The random variable $X(a,0)$ shall correspond to
the age of vertex $a$ at the time its first child $a1$ is born.
For $j\geq 1$, $X(a,j)$ shall parameterize the time between the
births of the $j$-th and $(j+1)$-th children of the $a$.
Therefore, the sequence $\{X(a,j)\}_{a,j}$ plays a role that is
similar to that of the $X_i$'s and $Y_j$'s in \subref{BinB} above.
There is, however, one important difference: whereas balls-in-bins
processes always have a fixed number of bins at which balls/hits
arrive, the number of ``bins" in the $\KR$ process grows. That is,
the potential vertices $a\in\natstr$ of the trees
$\{\tree_m\}_{m\geq 0}$ do not all come into existence at the same
time; they are rather {\em born} at appropriate times. We
therefore introduce a notion of {\em birth time}, which is defined
recursively as follows.
\begin{itemize} \item the birth time of the empty string $a=\estr$
is $\B(\estr)=0$; \item let $a_1, a_2,\dots,a_n\in\nat$ and
consider the sequence $a=a_1\dots a_{n}$. The birth time of $a$ is
the birth time of the parent sequence $b=a_1\dots a_{n-1}$ plus
the time until the $a_{n}$-th birth at $b$. More precisely,
$$\B(a)=\B(a_1\dots a_n) = \B(a_1\dots a_{n-1}) +
\sum_{j=0}^{a_n-1}X(a_1\dots a_{n-1},j)\periodeq$$
\end{itemize} An equivalent form of the definition of $\B(a)$ is
\begin{equation}\label{eq:KRintro_birth}\B(a)=\B(a_1\dots a_n) = \sum_{i=0}^{n-1}\sum_{j=0}^{a_{i+1}-1}X(a_1\dots a_{i},j)\periodeq\end{equation} Our
continuous time process is defined by setting $$\W(t)\equiv
\{a\in\natstr \; : \; \B(a)\leq t\} \;\;\;(t\in\real)\periodeq$$
$\W(\cdot)$ always takes values in the set $\E$ of~\proper~
subsets of $\natstr$ (defined in \subref{labelsmil}). This is
because the definition of birth time implies that the birth time
of $a_1\dots a_{n-1}$ is always
smaller than or equal to that of $a_1\dots a_n$.\\

Let us now specialize to the case where $f:\nat\cup\{0\}\to\real$
is given by $f(x)=(x+1)^p$ for some constant $p>1$. Such
attachment kernels satisfy the {\em explosion condition}
\begin{equation}\label{eq:explode}\sum_{n\geq 0}\frac{1}{f(n)}<+\infty\periodeq\end{equation}
The condition implies that the expectation of
\begin{equation}\label{eq:expl}\p(a) \equiv \sum_{j=0}^{+\infty}X(a,j) = \sup_{k\in\nat}\B(ak)-\B(a) \;\;\; (a\in\natstr)\end{equation}
is finite. Therefore, all the random variables defined in
\eqnref{expl} are almost surely finite.
\begin{definition} For an element $a\in\natstr$, the random variable $\p(a)$ defined in \eqnref{expl} is the {\em explosion time of $a$}. The infimum of $\B(a)+\p(a)$ over all $a\in\natstr$ is the {\em tree explosion time}, or the explosion time of the $\W(\cdot)$ process, and is denoted by $\s$.
\begin{equation}\label{eq:texpl}\s \equiv \inf_{a\in\nat}\B(a) + \p(a)\end{equation}\end{definition}

The intuition behind the definition of $\s$ is that it is the
first time when some node in the $\W(\cdot)$ process has an
infinite number of children. In fact, we {\em claim} that
\begin{claim}\label{cl:explodesright}\label{claim:explodesright}The following events hold with probability 1. The birth times $\B(a)$ that are smaller than $\s$ are pairwise distinct and can be well ordered with order type $\omega$.  Letting $0 = \B(\eps)
= \B_0 < \B_1 < \ldots < \B_n < \ldots $denote their ordered
sequence, $\B_n \nearrow \s$ as $n\to +\infty$.  Moreover, there
exists a unique $v \in \natstr$ that has infinite degree in
$\W(\s)$; this $v$ satisfies $\B(v)+\p(v)=\s$ and $\B(w)+\p(w)>
\s$ for all $w\neq v$.\end{claim}

A direct consequence of \clref{explodesright} is \thmref{expemb}
below.

\begin{theorem}\label{thm:expemb}Let $\{\B_n\}_{n\geq 0}$ be as in \clref{explodesright}. Then the sequence $\{\tree_n \equiv \W(\B_n)\}_{n\geq 0}$ is identical in law to the labelled $\KR$ process. Moreover, $\tree_n\to \ftree \equiv \W(\s)$ as $n\to +\infty$.\end{theorem}

However, to prove \clref{explodesright}, we will need some
elements of the proof of \thmref{KRintro_main1}. This could
potentially result in a problem: using \thmref{KRintro_main1} to
prove \clref{explodesright}, then employing the Claim to prove
\thmref{expemb}, and finally using this Theorem in the proof of
\thmref{KRintro_main1} would not be acceptable. Instead, we
circumvent this difficulty as follows.
\begin{enumerate}
\item In the beginning of the next Section, we state
\lemref{KRmain1_main1}, which is the same as
\thmref{KRintro_main1} but with $\W(\s)$ replacing $\ftree$ in the
statement. \item \clref{explodesright} is then proven, assuming
the Lemma.\item The remainder of the Section proves
\lemref{KRmain1_main1}, {\em without assuming}
\clref{explodesright} or \thmref{expemb} in any way. \item The
argument below shows that \clref{explodesright} implies
\thmref{expemb}, which directly implies that
\lemref{KRmain1_main1} can be strengthened to
\thmref{KRintro_main1}.
\end{enumerate}

Irrespective of formal proofs, the reader should keep in mind that
$\W(\s)$ represents the tree $\ftree$ in the statements of
\thmref{KRintro_main1} and \thmref{KRintro_main2}. Vertices
$a\in\natstr$ whose birth times satisfy $\B(a)>\s$ are not really
``born'' in $\ftree$, but rather constitute a {\em fictitious
continuation} of $\ftree$ in which new vertices continue to arrive
even though infinitely many vertices have already appeared. We
will use this continuation  to our advantage in many of the proofs
below.\\

\begin{proof}{[of \thmref{expemb}]} Assuming \clref{explodesright}, it suffices to show that for all $A\in\Efin$ and all $a\in A$
\begin{multline}\label{eq:pra}\Pr{\mbox{first birth of $\W(\cdot)$
after time $t$ is at $a$}\mid \W(t)=A} \\ =
\Pr{\B(a\deg_A(a)) + X(a,\deg_A(a)) = \min_{b\in A}\B(b\deg_A(b))+X(b,\deg_A(b))\mid \W(t)=A} \\
= \Pi(A,A\cup\{a(\deg_A(a)+1)\})\periodeq\end{multline} To prove
this, we first observe that the conditioning event is
\begin{equation*}\{\W(t) = A\} = \left\{\forall c\in A\;\; \B(c) + \sum_{j=0}^{\deg_A(c)-1}X(c,j)\leq t < \B(c) +
\sum_{j=0}^{\deg_A(c)}X(c,j)\right\}\periodeq\end{equation*} We
proceed as in the previous section and condition on the values
$X(b,j)=x(b,j)\geq 0$ for $b\in A$ and $0\leq j\leq \deg_A(b)-1$.
We want this event to be a subset of $\{W(t)=A\}$, so we require
that the birth times of all $b\in A$ are at most $t$; that is, we
must have:
\begin{equation}\label{eq:compcond}\forall b_1\dots b_r\in A\;\; y(b_1\dots b_r)\equiv
\sum_{i=0}^{r-1}\sum_{j=0}^{b_{i+1}-1}x(b_1\dots b_{i},j)\leq
t\periodeq\end{equation} Under this more stringent conditioning,
the probability we wish to compute is
\begin{equation}\label{eq:expprob}\Pr{\left. y(a\deg_A(a)) +
X(a,\deg_A(a)) =\min_{b\in A}y(b\deg_A(b))+X(b,\deg_A(b))\right|
B}\commaeq\end{equation} where $$B\equiv \{\forall b\in A\;\;
X(b,\deg_A(b))>t-y(b,\deg_A(b))\}\periodeq$$ The exponential
random variables in \eqnref{expprob} are all independent.
Moreover, by the lack of memory property, $X(b,\deg_A(b)) +
y(b,\deg_A(b)) - t$ conditioned on $X(b,\deg_A(b)) > t -
y(b,\deg_A(b))$ is distributed as $\exp(f(\deg_A(b)))$. It follows
that
\begin{multline}\label{eq:expprob2}\Pr{y(a\deg_A(a)) + X(a,\deg_A(a)) =\min_{b\in A}y(b\deg_A(b))+X(b,\deg_A(b))\mid B} \\
= \Pr{\exp\bigl(f(\deg_A(a))\bigr)-t = \min_{b\in
A}\exp\bigl(f(\deg_A(b))\bigr)-t}\end{multline}  where all
$\exp$'s are independent. From the minimum property, this last
probability is
\begin{equation}\frac{f(\deg_A(a))}{\sum_{b\in
A}f(\deg_A(b))}\end{equation} and this holds irrespective of the
values $\{x(b,j)\}\commaeq$ as long as \eqnref{compcond} is
satisfied. As a result, \eqnref{pra} holds. \end{proof}

\begin{remark}\label{rem:CTMC} The proof of \thmref{expemb} makes it clear that the $\W(\cdot)$ process is a continuous-time Markov Chain on $\Efin$ up to time $\s$. A consequence of this is the following. Let $t\geq 0$ be given and let $E$ be an event for $\W(\cdot)$ that is entirely defined in terms of $\{\W(s)\}_{0\leq s\leq t}$. Assume that inside the event $E$, $\W(t)=A\in\Efin$
. Finally, let $F$ be an event defined entirely in terms of
$\{\W(s)\}_{s\geq t}$. Then
$$\Pr{F\mid E} =\Pr{F\mid \W(t) = A}$$ We will employ this Remark in the proof of \lemref{finchild} and \lemref{posproblim} below.\end{remark}

\subsection{Three useful lemmas}\label{sub:lemmas}
Before we move on to prove the main theorems in the paper, we
collect three lemmas (proven in the Appendix) that will be useful
in dealing with sums of independent exponential random variables.
The present lemmas provide estimates of several probabilities that
are intimately related with the presence of nodes with $k$
descendants in the final tree $\ftree$. All of them are key
ingredients of the proofs of \thmref{KRintro_main1} and
\thmref{KRintro_main2}. We assume that $f(x)=(x+1)^p$ with $p>1$
in all statements.

\begin{lemma}[A large-deviations bound]\label{lem:lgdev}There exist constants $C,n_0>0$ depending only on $f$ such that for all $n\geq n_0$, all independent sequences of random variables $\{X_j\eqdist \exp(f(j))\}_{j\geq n}$,
and all $\delta>0$
\begin{eqnarray*}\Pr{\sum_{j\geq  n} X_j> \Ex{\sum_{j\geq n} X_j} + \delta} &\leq& C e^{- \delta n^{p-\frac 1 2}}\\
\Pr{\sum_{j\geq n} X_j< \Ex{\sum_{j\geq  n} X_j}  - \delta} &\leq&
C e^{- \delta n^{p-\frac 1 2}}\end{eqnarray*}
\end{lemma}

\begin{lemma}\label{lem:prprod}Let $Y=Y_1+\dots +Y_k$ be a sum of $k$ independent random variables for which
$$\forall 1\leq i\leq k\;\; \Pr{Y_i\leq \eps}=\Theta(\eps) \text{ as }\eps\searrow 0$$ Then there exist constants $C,n_0$ depending only on $k$, $f$ and the distributions $Y_i$ such that for all independent sequences of random variables $\{X_j\eqdist\exp(f(j))\}_{j\geq n}$ that are independent of $Y$ and all $n\geq n_0$
$$\frac{1}{C n^{k(p-1)}} \leq \Pr{Y\leq \sum_{j>n} X_j }\leq \frac{C}{n^{k(p-1)}}$$
\end{lemma}
\begin{lemma}\label{lem:Poi} Let $Z_1,\ldots,Z_k$ be
independent exponentials with mean one and let $Z\equiv
Z_1+\ldots+Z_k$. Then for all $\lambda>0$
$$\Pr{Z\leq \lambda} = e^{-\lambda}\sum_{j\geq k} \frac{\lambda^j}{j!}\leq \frac{\lambda^k}{k!}$$\end{lemma}

\section{Finitely many k-fertile vertices}\label{sec:KRmain1}
In this section we prove the first of our main results about the
$\KR$ process, \thmref{KRintro_main1}. As noted in the previous
section, \clref{explodesright} -- which has not been proven yet --
is necessary for the connection between the exponential process
and the $\KR$ process. Proving the Claim will require the a
preliminary form of \thmref{KRintro_main1} that we shall present
below. We assume throughout the section that $f(x)=(x+1)^p$ for
some $p>p_k=1+1/k$.\\

Recall that a node is {\em $k$-fertile} if it has $k$ or more
descendants in the corresponding tree.

\begin{lemma}\label{lem:KRmain1_main1}Consider the $\W(\cdot)$ process defined in \secref{KRintro_expemb}, and assume its attachment kernel is $f(x)=(x+1)^p$, $p>p_k$. Then, for all $T> 0$,
\begin{equation}\label{eq:KRmain1_finalaim}\Ex{\#\{a\in\W(\s): a\mbox{ is $k$-fertile}\}\mid \p(\estr)\leq
T}<+\infty\commaeq\end{equation} and therefore
$$\Pr{\#\{a\in\W(\s): a\mbox{ is $k$-fertile}\}<+\infty\mid
\p(\estr)\leq T}=1\periodeq$$ Since $\p(\estr)<+\infty$ almost
surely, this implies that $\W(\s)$ almost surely has only finitely
many $k$-fertile vertices.\end{lemma}

As noted in \secref{KRintro_expemb}, we will use \lemref{KRmain1_main1} to prove \clref{explodesright}, and this in turn will imply that \thmref{expemb} holds. This last Theorem and \lemref{KRmain1_main1} directly imply \thmref{KRintro_main1}. Therefore, most of the present section will be devoted to proving \lemref{KRmain1_main1}.\\

This section is organized as follows. In \subref{KRmain1_claim} we
show how \lemref{KRmain1_main1} implies \claimref{explodesright}.
Having settled that matter, we move on to proving
\lemref{KRmain1_main1}. Our proof will consist of bounding the
probabilities of the form $$\Pr{a\mbox{ is $k$-fertile}\mid
\p(\estr)\leq T}\commaeq$$ and then showing that their sum is
finite. We illustrate our techniques for doing so in
\subref{KRmain1_easybound} below, where we show a partial result
in the direction of \lemref{KRmain1_main1}. We then show in
\subref{KRmain1_subtree} that the time at which a given
$a\in\Nstr$ becomes $k$-fertile in the $\W(\cdot)$ process can be
bounded in terms of a sum of $k$ exponential random variables
(\lemref{KRmain1_kdesc}). This permits an improved bound on the
probability of $k$-fertility (\subref{KRmain1_mixed}), which is
then applied to prove \lemref{KRmain1_main1} in
\subref{KRmain1_final}.

\subsection{\lemref{KRmain1_main1} implies \claimref{explodesright}}\label{sub:KRmain1_claim}

\begin{proof}{[of \claimref{explodesright}]} The following lemma is a well-known combinatorial result.

\begin{lemma}{[K\" onig's Infinity Lemma]}\label{lem:KRmain1_konig} Let $T$ be an infinite rooted tree in which every vertex has finite degree. Then $T$ contains an infinite path starting from the root.\end{lemma}

We will use the Infinity Lemma and \lemref{KRmain1_main1} to prove
a series of almost-sure statements that imply the Claim.\\

{\em All birth times are almost surely distinct.} This occurs
because, for all distinct $a,b\in\Nstr$, the difference
$\B(a)-\B(b)$ is a sum of terms of the form $\pm X(c,j)$ for some
$(c,j)\in\Nstr\times \N\cup\{0\}$. Each such term has a smooth
distribution with no point masses, and all terms are independent,
hence $\B(a)-\B(b)\neq 0$ with probability $1$.\\

{\em There almost surely exists at least one vertex $v\in\Nstr$
with infinite degree in $\W(\s)$.}For suppose that this were note
the case. Since $\W(\s)$ is infinite, the Infinity Lemma would
imply that there was an infinite path starting from the root in
$\W(\s)$. But all the infinitely many vertices on such path would
have $\geq k$ descendants, for any $k\in\N$. However, $p>1$
implies that $p>p_k = 1+1/k$ for {\em some} $k\in\N$, and
\lemref{KRmain1_main1} then implies that only finitely many
vertices in $\W(\s)$ can be
$k$-fertile, a contradiction.\\

{\em There almost surely exists a unique vertex $v$ for which
$\s=\B(v) + \p(v)$.} With probability $1$, there is a vertex $v$
of $\W(\s)$ with infinite degree. Since the degree of $v$ is
infinite in $\W(\s)$, all the children of $v$ must have been born
before time $\s$.
$$\forall n\in\N\commaeq\;\; \B(vn)\leq \s\periodeq$$
As $n\to +\infty$, $\B(vn)\to \B(v)+\p(v)$, thus $\B(v)+\p(v)\leq
\s$. Then, by definition of $\s$, $\s=\B(v)+\p(v)$. Thus there
{\em exists} a $v$ as claimed. For uniqueness, one can show that
$\B(a)+\p(a)\neq \B(b)+\p(b)$ for all distinct $a,b\in\Nstr$.\\

{\em With probability $1$, $\W(t)$ is finite for all $t<s$.}
Suppose that is not the case. For all $a\in \Nstr$
$$\B(a)+\p(a) =\lim_{d\to +\infty}\B(ad) \geq \s>t\commaeq$$
which implies that for all $a$ there is an integer $d_a\geq 0$
such that $\B(ad_a)>t$. Therefore, any $a$ has finite degree $\leq
d_a-1$ in $\W(t)$. By the Infinity Lemma, $\W(t)$ must then have
an infinite path from $\estr,a_1,a_1a_2,a_1a_2a_3,\dots$. But all
nodes along this path have infinitely many descendants in $\W(t)$,
and hence also in $\W(\s)$, which was shown above to have probability $0$. The contradiction implies the assertion.\\

{\em The set of birth times before $\s$ can be well-ordered.} This
is a consequence of the previous assertion.\\

{\em With probability $1$, when the descendants of $v$ are removed
from $\W(\s)$, the result is a finite tree.} Again, the key
property here is that all $a\in\Nstr\backslash\{v\}$ have finite
degree. So if $\W(\s)$ without the descendants of $v$ would be
infinite, the Infinity Lemma would imply the existence of an
infinite path in $\W(\s)$, which would imply that all nodes along
the path have infinitely many descendants. Since this is
impossible, the assertion must be
true.\\

{\em The ordered birth times $\B_0=0\leq \B_1\leq \B_2\leq \dots$
are almost surely distinct and converge almost surely to $\s$.}
That they are distinct follows from the first assertion. Since
they form an increasing sequence bounded by $\s<+\infty$, they
converge to some finite limit. But the birth times
$\{\B(vm)\}_{m=1}^{+\infty}$ (with $v$ as in the previous
paragraph) form a subsequence of $\{\B_n\}_{n\in\N\cup\{0\}}$ that
converges to $\B(v)+\p(v)=\s$,
so the $\{\B_n\}_{n}$ sequence converges to $\s$ as well.\\

The series of assertions implies the Claim.\end{proof}

\subsection{Two instructive examples}\label{sub:KRmain1_easybound}

Having shown that \lemref{KRmain1_main1} implies
\claimref{explodesright}, we now turn to the proof of the Lemma.
Recall that the goal of that lemma is to prove that only finitely
many vertices have $k$ or more descendants in $\W(\s)$. For the
sake of the reader, however, we first consider two special classes
of $a\in\Nstr$ and prove that only finitely many nodes in each
class have large degree. While the corresponding general result
combines ingredients of the two special cases below, we believe
that our techniques become much clearer if introduced separately.\\

To state the present results, we need two definitions. Fix a
number $L>0$, and call $a\in\Nstr\backslash\{\estr\}$ {\em
$L$-moderate} if all numbers in the sequence $a$ are smaller than
or equal to $L$. If on the other hand all numbers in $a$ are
bigger than $L$, call it {\em $L$-extreme}. Our two simple lemmas
are presented below.

\begin{lemma}\label{lem:KRmain1_Lmod} For all integers $L>0$ and all $T\geq 0$, the expected number of $L$-moderate $1$-fertile vertices in $\W(\s)$ conditioned on $\p(\estr)=T$ is finite.\end{lemma}
\begin{lemma}\label{lem:KRmain1_Lext} There exists a constant $L_0>0$ defined only in terms of $p$ such that for all integers $L\geq L_0$, the expected number of $L$-extreme vertices in $\W(\s)$ that have at least $k$ children is finite.\end{lemma}

\begin{proof}{[of \lemref{KRmain1_Lmod}]} For any
$a=a_1\ldots a_m$, the time for the birth of the first child of
$a$ is \begin{multline}\B(a1)= \B(a) + X(a,0) \\ =
\sum_{j=0}^{m-1}\sum_{i=0}^{a_{j+1}-1} X(a_1\ldots a_j,i) +
X(a_1\dots a_m,0)\geq \sum_{j=1}^m X(a_1\ldots
a_j,0)\periodeq\end{multline} Notice that this lower bound on
$\B(a1)$ is actually {\em independent} of $\p(\estr)$, which is at
least as big as the tree explosion time $\s$. As a result:
\begin{multline*}\Pr{a\mbox{ is $1$-fertile}\mid \p(\estr)=T} =
\Pr{\B(a1)\leq \s
\mid \p(\estr)=T} \\
\leq \Pr{\sum_{j=1}^m X(a_1\dots a_j,0)\leq \p(\estr)\mid
\p(\estr)=T} = \Pr{\sum_{j=1}^m X(a_1\ldots a_j,0)\leq
T}\periodeq\end{multline*} We now apply \lemref{Poi} with $\lambda
= T$ and $Z = \sum_{j=1}^m X(a_1\ldots a_j,0)$ to deduce
$$\Pr{a \mbox{ is $1$-fertile }\mid \p(\estr)=T}\leq \Pr{\sum_{i=1}^m
X(a_1\ldots a_i,1)\leq T} \leq \frac{T^m}{m!}$$ There are $L^m$
$L$-moderate $a$ of length $m$, and this implies that
\begin{equation*}\sum\limits_{a\mbox{ $L$-moderate}}\Pr{a\mbox{ is
$1$-fertile} \mid\p(\estr)=T} \leq\sum_{m=1}^{+\infty}
\frac{T^mL^m}{m!} = e^{TL}-1<+\infty\periodeq\end{equation*} This
finishes the proof.\end{proof}

\begin{proof}{[of \lemref{KRmain1_Lext}]} We assume $L\geq k+1,n_0$, where $n_0$ comes from \lemref{prprod}. Fix an $L$-extreme $a=a_1\dots a_m$ with all $a_i\in\N$, and let $a_{m+1}=k$. The
event \begin{equation}\{a\mbox{ has at least $k$ children}\} =
\{\B(ak)\leq \s\} \end{equation} is contained the event
\begin{equation}H_{a,i}\equiv \left\{\sum_{j=0}^{k-1}X(a_1\dots a_i,j)\leq \sum_{s=a_{i}}^{+\infty}X(a_1\dots a_{i-1},s)\right\}\commaeq\end{equation}
for each $1\leq i\leq m$. This is true because for all $1\leq
i\leq m$
\begin{eqnarray*}\B(a_1\dots a_{i}) + \sum_{j=0}^{k-1} X(a_1\dots a_{i},j) &=& \B(a_1\dots a_{i}k) \\ & & \text{ (by properties of birth times)}\\
&\leq &\B(a_1\dots a_{i}a_{i+1}) \\
& & \text{   (since $a_i \geq L \geq k$ for $i<m$, and $a_{m+1}=k$)}\\
&\leq & \B(ak)\\
& &\text{ (since $ak$ is either a descendant of } \\
& &\text{ $a_1\dots a_{i+1}$ or equal to $ak$)}\end{eqnarray*} and
\begin{eqnarray*}\s  & \leq & \B(a_1\dots a_{i-1}) + \p(a_1 \dots a_{i-1}) \\ & & \text{ (by definition of $\s$)}\\
& = & \B(a_1\dots a_{i-1}) + \sum_{j=0}^{+\infty}X(a_1\dots a_{i-1},j) \\
& = & \B(a_1\dots a_{i}) + \sum_{j=a_{i}}^{+\infty}X(a_1\dots
a_{i-1},j)\periodeq\end{eqnarray*} So that
\begin{eqnarray*} & & \{\text{a has at least $k$ children}\}\\ & =&\{\B(ak)\leq \s \}\\
& \subseteq & \left\{\B(a_1\dots a_{i}) + \sum_{j=0}^{k-1} X(a_1\dots a_{i},j) \leq \B(a_1\dots a_{i}) + \sum_{j=a_{i}}^{+\infty}X(a_1\dots a_{i-1},j)\right\} \\
& = &\left\{\sum_{j=0}^{k-1} X(a_1\dots a_{i},j) \leq
\sum_{j=a_{i}}^{+\infty}X(a_1\dots a_{i-1},j)\right\} =
H_{a,i}\periodeq\end{eqnarray*} Now note that all the events
$\{H_{a,i}\}_{1\leq i\leq m}$ are in fact {\em independent}. In
fact, for any $1\leq i\leq n-1$, $H_{a,i}$ depends only on the
random variables $X(a_1\dots a_{i-1},j)$ with $j\geq a_i\geq L\geq
k+1$ and $X(a_1\dots a_i, \ell)$ with $0\leq \ell \leq k$.
Therefore, the choice of $L$ implies that no random variable can
appear in the definitions of two different $H_{a,i}$. Therefore,
\begin{equation}\nonumber\Pr{a\mbox{ has at least $k$ children}}\leq \prod_{i=1}^m\Pr{H_{a,i}}\periodeq\end{equation} Now notice that
\begin{equation}\Pr{H_{a,i}} = \Pr{\sum_{\ell=1}^{k}Y_\ell \leq \sum_{j\geq a_i}X(a_1\dots a_{i-1},j)}\commaeq\end{equation}
with $Y_\ell = X(a_1\dots a_i,\ell-1)$. It is straightforward to
check that the assumptions of \lemref{prprod} hold (since we know
$a_i\geq L\geq n_0$) and that as a
result\begin{equation}\label{eq:xxxxx}\Pr{H_{a,i}}\leq
\frac{C}{a_i^{(p-1)k}}\commaeq\end{equation}  where $C$ depends
only on $p$, as the distributions of the $Y_\ell$'s are determined
by $p$. It follows that
\begin{equation}\sum\limits_{a\mbox{ $L$-ext.}} \Pr{a\mbox{ has
$k$ children}}\leq \sum_{m=1}^{+\infty} \sum_{a_1,\dots,a_m
>L}\prod_{i=1}^m \frac{C}{a_i^{k(p-1)}} = \sum_{m=1}^{+\infty}
\left(\frac{C}{L^{k(p-1)-1}}\right)^m\periodeq\end{equation}
Noting that $p>p_k\Rightarrow k(p-1)>1$, we can now take $L \geq
L_0\equiv {(2C)}^{\frac{1}{k(p-1)-1}}$ to have a finite
sum.\end{proof}

\begin{remark}\label{rem:probkchildren} One can show by the same proof technique as above, that for all fixed $v\in\Nstr$ and all fixed $k\in\N$
\begin{equation}\label{eq:remprob}\Pr{vi\mbox{ has $k$ children before $v$ explodes} } = \bigoh{i^{-(p-1)k}}\mbox{ as }i\to +\infty\periodeq\end{equation}
To prove this, note that the event in \eqnref{remprob} is
\begin{equation*}\{\B(vik)\leq \B(v)+\p(v)\} = \left\{\sum_{j=0}^{k-1}X(vi,k)\leq \sum_{j\geq i}X(v,j)\right\}\commaeq\end{equation*}
because $\B(vik)=\B(vi) + \sum_{j=0}^{k-1}X(vi,k)$ and
$\B(v)+\p(v)=\B(vi)+\sum_{j\geq i}X(v,j)$. Then apply
\lemref{prprod}, as in the previous proof.\\

Similarly, one can show that, for all $t\geq 0$, all $v,w\in
\Nstr$ and all finite trees $\tree_n$ such that $v$ has $n$
children in $\tree_n$,
\begin{multline*}\Pr{w\mbox{ has $k$ children after time $t$ and before $v$ explodes}\mid \W(t)=\tree_n } \\ = \bigoh{n^{-(p-1)k}}\mbox{ as }n\to +\infty\periodeq\end{multline*}
We will employ this remark in the proof of
\lemref{posproblim}.\end{remark}

\subsection{Subtrees and the time until $k$ descendants are born}\label{sub:KRmain1_subtree}

There are two reasons why \lemref{KRmain1_Lmod} and
\lemref{KRmain1_Lext} do not imply \lemref{KRmain1_main1}. First,
there are $a\in\Nstr$ that are neither $L$-moderate nor
$L$-extreme. Second, the above lemmas only bound the probability
of a certain node having degree $\geq k$, which is different from
$k$-fertility for all $k\geq 2$. The next Lemma deals with the
latter difficulty. Fix some $a\in\Nstr$ and let $\W_a(t)\equiv
\{c\in\Nstr \; : \; \B(ac)-\B(a)\leq t\}$ (for $t\in\R$) be the
subtree of $\W(t+\B(a))$ rooted at $a$. Clearly, $\W_a(\cdot)$ and
$\W(\cdot) = \W_\estr(\cdot)$ have the same distribution.
Moreover, $a$ is $k$-fertile if and only if the size of
$\W_a(\s-\B(a))$ is at least $k+1$ (i.e. $\W_a(\s-\B(a))$ has at
least $k$ vertices other than the root). \lemref{KRmain1_kdesc}
provides tools for the analysis of the $k$-fertility event.

\begin{lemma}\label{lem:KRmain1_kdesc}For a fixed $a\in\Nstr$, let $\T_0(a)$ be the time of the first birth of a node other than the root in the $\W_a(\cdot)$ process. Moreover, for $i\in\N$, let $\T_i(a)$ be the time elapsed between the $i$th and $(i+1)$th births in $\W_a(\cdot)$ (again excluding the birth time of the root). Then there exist a sequence of random variables
$\{\sR_j(a)\}_{j=0}^{+\infty}$ such that:
\begin{enumerate}
\item $\{\sR_j(a)\}_{j=0}^{+\infty}$ is a sequence of independent
random variables; \item the sequence
$\{\sR_j(a)\}_{j=0}^{+\infty}$  is a deterministic function of the
random variables $\{X(ac,i)\mid c\in\Nstr,\,i\in\N\}$; \item for
each $j\in\N\cup\{0\}$, $\sR_j(a)\eqdist \exp((j+1)f(j))$; \item
$\sR_0(a)=\T_0(a)$ and for all $j\in\N$ $\sR_j(a)\leq \T_j(a)$.
\end{enumerate}\end{lemma}
\begin{proof}{} It suffices to consider the case $a=\estr$. For convenience, we introduce the notation
$$\Sigma(c)\equiv \left\{\begin{array}{ll}\sum_{i=1}^{m}c_i & c=c_1\dots c_m \in\Nstr\backslash\{\estr\}\\
0 & c = \estr\end{array}\right.\periodeq$$ We prove inductively
that the random variables $\{\sR_j(\estr)\}_{j=0}^{r}$ can be
defined as above, so that for all $j\in\N\cup\{0\}$ $\sR_j(\estr)$
is completely defined by the values of $X(c,r)$ for $c\in\Nstr$,
$0\leq \Sigma(c)+j\leq r$. For $r=0$, this is easy: just set
$\sR_0(\estr)=\T_1(\estr)=X(\estr,0)$. Now assume inductively that
$\sR_j(\estr)$ has been defined for all $0\leq j\leq r=n-1$. To
prove that the same is possible for $r=n$, condition on a
particular value
\begin{equation}\label{eq:fertimes_cond}\W_\estr\left(\sum_{j=0}^{n-1}\T_j(\estr)\right)=A\in\E\commaeq |A|<+\infty\periodeq\end{equation}
$\sum_{j=1}^{n-1}\T_j(\estr)$ is exactly the birth time of the
$n$th descendant of the root in $\W_\estr(\cdot)$ (for $\estr$ is
born at time $0$), hence $|A|=n+1$. We also notice that,
$\Sigma(c)+\deg_A(c)\leq n$ for all $c\in A$. Indeed, the sequence
$$b=\left\{\begin{array}{ll}c\deg_A(c) & \mbox{ if }\deg_A(c)>0 \\ c & \mbox{ if }\deg_A(c)=0\end{array}\right.$$ is an element of $A$ with $\Sigma(b)=\Sigma(c)+\deg_A(c)$, and it is a simple fact (whose proof we omit) that $\Sigma(b)\leq |A|-1$ for {\em any} $b\in A\in\Efin$.\\

Conditioned on the event in \eqnref{fertimes_cond}, the random
variable $\T_{n}(a)$ has exponential distribution with rate
$\sum_{c\in A}f(\deg_A(c))$, which is bounded by $|A|f(n)\leq
(n+1)f(n)$ by the above remarks. Therefore,
\begin{equation*}\sR_{n}(a)\equiv \frac{\sum_{c\in A}f(\deg_A(c))}{(n+1)f(n)}\T_{n}(a) \leq \T_{n+1}(a) \;\;
\left(\mbox{where }A =
\W_\estr\left(\sum_{j=1}^{n-1}\T_j(\estr)\right)\right)\end{equation*}
is exponential with rate $(n+1)f(n)$ irrespective of $A$, by the
multiplication property of exponentials (cf. \subref{probing}).
Because $\W_a(\sum_{j=1}^{n-1}\T_j(a))$ and $\T_{n}(a)$ are
completely defined by the random variables $\{X(ac,j)\; : \;
c\in\Nstr, \Sigma(c)+j\leq n+1\}$, the same is true of
$\sR_{n}(a)$. This finishes the proof.\end{proof}

\subsection{A general bound on the probability of $k$-fertility}\label{sub:KRmain1_mixed}

\lemref{KRmain1_kdesc} is now used to prove a stronger form of the
bounds in \subref{KRmain1_easybound} that applies to {\em all}
$a\in\Nstr$ (and not just $L$-moderate or $L$-large sequences). To
present this bound, we need a definition. For a fixed $L>0$ and a
sequence $a=a_1a_2\dots a_m\in\Nstr$ of length $m$, the set of
{\em small indices in $a$} is $\smL(a)=\{1\leq i\leq m\; : \;
a_i\leq L\}$, and the set of {\em large indices in $a$} is
$\lgL(a)=\{1\leq i\leq m\; : \; a_i> L\}$.
\begin{lemma}\label{lem:KRmain1_mixed}There exist constants $C,L_0>0$ depending only on $k$ and $p$ such that for any $T>0$, $L\geq L_0$ and $a=a_1\dots a_m\in\Nstr$ \begin{equation}\label{eq:mixed}\Pr{a\mbox{ is
$k$-fertile and } \p(\estr)\leq T} \leq
\frac{\max\{T,C\}^{m}}{|\smL(a)|!} \,\prod_{j\in\lgL(a)}
\frac{1}{a_j^{(p-1)k}}\periodeq\end{equation}\end{lemma}
\begin{proof}{} For most of the proof, we will only assume that $L_0\geq k$; more conditions on $L_0$ will be imposed later. Set $a_{m+1}\equiv k$ and for each $i\in\lgL(a)$ define $\hati$ to be the smallest
$j\in\lgL(a)\cup\{m+1\}$ satisfying $j>i$; notice that the choice
of $L_0$ implies $a_{\hati}>L_0\geq k$ whenever $\hati<m+1$ .
Employing the random variables $\{\sR_j(a)\}_{j=0}^{k-1}$ whose
existence \lemref{KRmain1_kdesc} guarantees, we deduce that
\begin{eqnarray}\{a\mbox{ is $k$-fertile}\} &\subset & \{\B(a) +
\sum_{n=0}^{k-1}\sR_n(a)\leq \s\} \\ \nonumber &=& \{\B(a) +
X(a,0) + \sum_{n=1}^{k-1}\sR_n(a)\leq
\s\}\periodeq\end{eqnarray}In what follows, we will bound the
probability on the right-hand side, noting that the $\sR_j(a)$'s
and $X(b,i)$'s that appear in the definitions below are all
independent because of \lemref{KRmain1_kdesc}. Consider the
following events.

\begin{eqnarray}\label{eq:kfert} F^T_a & \equiv & \left\{ \B(a) + X(a,0) + \sum_{n=1}^{k-1}\sR_n(a)\leq \s\mbox{ and } \p(\estr)\leq T\right\}\commaeq\end{eqnarray}
\begin{eqnarray}\label{eq:kfert11} G^T_a& \equiv & \left\{\sum_{i\in\smL(a)}X(a_1\dots a_i,0) \leq
T\right\}\commaeq\end{eqnarray}
\begin{equation}H_{a,i} \equiv  \left\{\left(\begin{array}{l} X(a_1\dots
a_i,0) \\ + \sum_{j=1}^{k-1} X(a_1\dots a_{\hati
-1},j)\end{array}\right)\leq \sum_{j\geq a_i}X(a_1\dots
a_{i-1},j)\right\} \end{equation}
\begin{equation*}\;\;(i\in\lgL(a), \hati\neq m+1) \commaeq\end{equation*}
\begin{equation}\label{eq:kfert33} H_{a,i} \equiv
\left\{X(a_1\dots a_i,0)+\sum_{j=1}^{k-1} \sR_j(a)\leq \sum_{j\geq
a_i}X(a_1\dots a_{i-1},j)\right\}\end{equation}
\begin{equation*}\;\;(i\in\lgL(a), \hati = m+1)\periodeq \end{equation*}
The first event is the one whose probability we want to bound. The
second event is similar to the one in the proof of
\lemref{KRmain1_Lmod}, whereas the remaining events are
reminiscent of those in the proof of \lemref{KRmain1_Lext}. We now
{\em claim} that:

\begin{claim}\label{claim:1}It holds that
\begin{equation}\label{eq:kfert2}F^T_a \subseteq G^T_a \cap \bigcap\limits_{i\in\lgL(a)}H_{a,i}\periodeq\end{equation}
Moreover, the events on the right-hand side of \eqnref{kfert2} are
independent.\end{claim}

\claimref{1} is proven at the end of the current proof. but we now
present the following concrete example of its application to
illustrate our argument. Assume that $k=2$, $L=L_0=3$ and
$a=a_1a_2\dots a_6=142461$, in which case $\smL(a)=\{1,3,6\}$ and
$\lgL(a)=\{2,4,5\}$. \figref{claim2} represents some of the random
variables involved in \eqnref{kfert2} by rectangles. The first six
columns of rectangles stand for random variables of the form
$X(b,j)$ for $b=\estr$ (the empty string), $1$, $14$, $\dots$,
$14246$, and $j=0,1,\dots,6$, while the last column represents the
random variables $\sR_0(a)=X(a,0)$ and $\sR_1(a)$. The rectangles
that lie completely below the dashed line correspond to the random
variables that appear in
\begin{figure}[t]\begin{center}
\includegraphics[scale=.80]{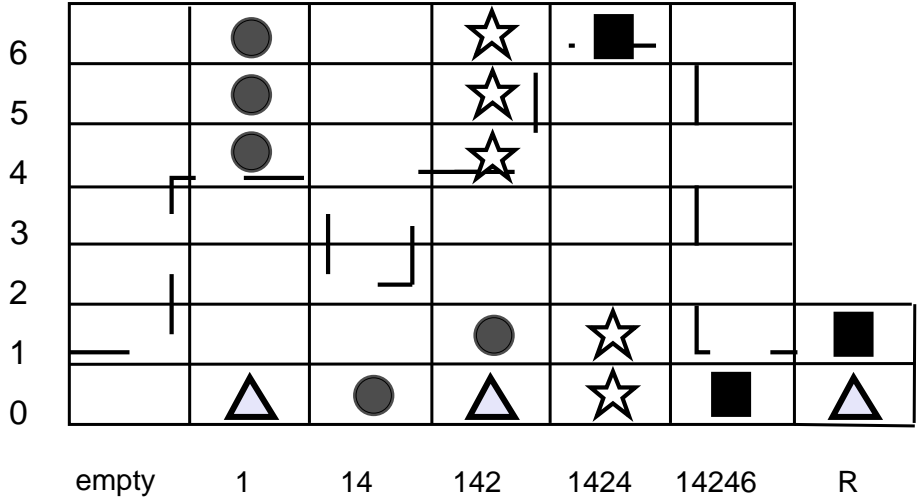}\end{center}\caption{\normalsize
Illustration of \claimref{1} for $a=142461$.}\label{fig:claim2}
\end{figure}
\begin{eqnarray}\label{eq:Bacontains}\B(a)  &=& X(\estr,0) + \sum_{j_1=0}^{3}X(1,j_1)
+ \sum_{j_2=0}^1 X(14,j_2) \\ \nonumber& & + \sum_{j_3=0}^3
X(142,j_3) + \sum_{j_4=0}^5 X(1424,j_4) +
 X(14246,0)\periodeq\end{eqnarray}
Moreover,
\begin{equation*}\hati =
\left\{\begin{array}{ll} 4\commaeq & i=2 \\ 5\commaeq & i=4 \\
7\commaeq & i=5\periodeq\end{array}\right.\end{equation*} By
checking the definitions of $G^T_a$ and $H_{a,i}$, one can check
that the following statements hold.
\begin{enumerate}
\item {\em $G^T_a \supseteq F^T_a$}, since in the event $F^T_a$
the explosion time $\s$ is at most $T$, and the sum
\eqnref{Bacontains} defining $\B(a)+X(a,0)$ contains the terms of
$\sum_{i\in\smL(a)}X(a_1,\dots a_i,0)$. Moreover, the random
variables appearing in $G^T_a$ correspond to the rectangles marked
with triangles in \figref{claim2}. \item {\em $H_{a,2}\supseteq
F^T_a$}. In order for $F^T_a$ to happen, $a_1a_2\dots
a_{I_2}=1424$ must be born before $a_1=1$ explodes. In particular,
 using the critical fact that $k=2<L=3$, so that $a_{I_2}>k$,
$1422=a_1a_2\dots a_{I_2-1}k$ must be born before node $a_1=1$
explodes . Since $a_1$ explodes at time
$\B(a_1)+\p(a_1)=\B(1)+X(1,0)+X(1,1)+\dots$, $1422$ is born at a
time that is larger than $\B(14)+X(14,0)+X(142,1)$ and $\B(14) =
\B(1) + X(1,0)+\dots +X(1,3)$, it follows that
\begin{multline*}F^T_a \\
\subset\left\{\left(\begin{array}{l}\B(1)+ X(1,0)+ \\
+\dots + X(1,3) + \\
X(14,0)+X(142,1)\end{array}\right) \leq \B(1)+X(1,0)+X(1,1)+\dots\right\} \\
\subset\left\{X(14,0)+X(142,1) \leq \sum_{j\geq 4} X(1,j)\right\}=
H_{a,2}\commaeq\end{multline*} so $H_{a,2}\supseteq F^T_a$ is
indeed true. Moreover, one can check that the random variables
appearing in the definition of $H_{a,2}$ are precisely the ones
marked with circles in \figref{claim2}. \item Similarly, one can
show that $H_{a,4}\supseteq F^T_a$ (respectively,
$H_{a,5}\supseteq F^T_a$) and that the random variables marked
with stars (resp. squares) are precisely the ones appearing in the
definition of $H_{a,4}$ (resp. $H_{a,5}$).
\end{enumerate}
Items $1.$, $2.$ and $3.$ above imply not only the validity of
\eqnref{kfert2}, but also that {\em no random variable of the form
$X(\cdot,\cdot\cdot)$ or $\sR_{\cdot}(a)$ appears in the
definition of more than one of the events in
\eqnref{kfert11}--\eqnref{kfert33}}. Since those random variables
are also {\em independent}, we have proven that $G^T_a$,
$H_{a,2}$, $H_{a,4}$ and $H_{a,5}$ are independent events, which
implies the Claim in this special case. The proof of \claimref{1}
for general $L$, $k$ and $a$ is entirely
analogous to the argument sketched above.\\

We continue with the proof of \lemref{KRmain1_mixed}, noting that
\claimref{1} implies
\begin{equation}\label{eq:boundind}\Pr{F^T_a}\leq \Pr{G^T_a} \times \prod_{i\in\lgL(a)} \Pr{H_{a,i}}\periodeq\end{equation}
The remainder of our proof consists of bounding the probabilities
on the right-hand side of \eqnref{boundind}, which is done in
roughly the same way as in Lemmas $\ref{lem:KRmain1_Lmod}$ and
$\ref{lem:KRmain1_Lext}$. The probability of $G^T_a$ is bounded
using \lemref{Poi} with the $Z_i's$ corresponding to the
$X(a_1\dots a_i,0)$ for $i\in\smL(a)$ and $\lambda=T$.
\begin{equation}\label{eq:boundG}\Pr{G^T_a}\leq \frac{T^{|\smL(a)|}}{|\smL(a)|!}\periodeq\end{equation}
Now fix some $i\in\lgL(a)$ with $\hati\neq m+1$. We apply
\lemref{prprod} with $Y_1=X(a_1\dots a_i,0)$, $Y_\ell = X(a_1\dots
a_{\hati-1},\ell-1)$ ($2\leq\ell \leq k$) and $$\{X_j\}_{j\geq n}
= \{X(a_1\dots a_{i-1},j)\}_{j\geq a_i}\periodeq$$ In the present
case, the distributions of the $Y_\ell$'s are all defined in terms
of $f$ and $k$. Therefore there exist $C,n_0$ depending only on
$k$ and $f$ such that if $a_i\geq n_0$,
\begin{equation}\label{eq:boundHf}\Pr{H_{a,i}}\leq \frac{C}{a_i^{(p-1)k}}\periodeq\end{equation}
For $i\in\lgL(a)$ with $\hati=m+1$, a similar reasoning with
$Y_\ell = \sR_\ell(a)$ for $2\leq \ell\leq k$ implies that for
(possibly enlarged) $C,n_0$ depending only on $k$ and $f$, and all
$a_i\geq n_0$, \eqnref{boundHf} still holds. So if we take $L_0
\geq n_0$, we can plug \eqnref{boundG} and \eqnref{boundHf} into
\eqnref{boundind} for any $a$, which finishes the
proof.\end{proof}

To conclude, we now prove \claimref{1}.\\

\begin{proof}{[of \claimref{1}]}We first show
that show that each of the events that are (re)defined below
\begin{eqnarray}\label{eq:KRmain1_Gagain}G^T_a& \equiv & \left\{\sum_{i\in\smL(a)}X(a_1\dots a_i,0) \leq
T\right\}\commaeq\end{eqnarray}
\begin{equation}\label{eq:KRmain1_Hagain}H_{a,i} \equiv  \left\{\left(\begin{array}{l} X(a_1\dots
a_i,0) \\ + \sum_{j=1}^{k-1} X(a_1\dots a_{\hati
-1},j)\end{array}\right)\leq \sum_{j\geq a_i}X(a_1\dots
a_{i-1},j)\right\} \end{equation}
\begin{equation*}\;\;(i\in\lgL(a), \hati\neq m+1)\commaeq \end{equation*}
\begin{equation}\label{eq:KRmain1_Hagain2}H_{a,i} \equiv
\left\{X(a_1\dots a_i,0)+\sum_{j=1}^{k-1} \sR_j(a)\leq \sum_{j\geq
a_i}X(a_1\dots a_{i-1},j)\right\}\end{equation}
\begin{equation*}\;\;(i\in\lgL(a), \hati = m+1)\commaeq \end{equation*}
contains
\begin{equation}\label{eq:defF22}F^T_a =\left\{ \B(a) + X(a,0) + \sum_{j=1}^{k-1}\sR_j(a)\leq \s\mbox{ and } \p(\estr)\leq T\right\}\periodeq\end{equation}
\mypar{First containment: $G^T_a\supset F^T_a$. }On the one hand,
all terms appearing in the sum $$\sum_{i\in\smL(a)}X(a_1\dots
a_i,0)$$ also appear in the sum defining $\B(a)=\B(a_1\dots a_m)$
(cf. \eqnref{KRintro_birth}), so that
\begin{equation}\sum_{i\in\smL(a)}X(a_1\dots a_i,0) \leq \B(a) \leq \B(a) + X(a,0) + \sum_{j=1}^{k-1}\sR_j(a)\periodeq\end{equation}
On the other hand, by the definition \eqnref{texpl} of $\s$
\begin{equation}\s \leq \B(\estr) + \p(\estr) = \p(\estr) \text{   since  }\B(\estr)=0\periodeq\end{equation}
Therefore,
$$F^T_a\text{ occurs } \Rightarrow \sum_{i\in\smL(a)}X(a_1\dots a_i,0)  \leq \p(\estr)\leq T \Rightarrow G^T_a \text{  occurs  }\periodeq$$

\mypar{Second containment: $H_{a,i}\supset F^T_a$ if $i\in\lgL(a)$
and $\hati<m+1$. }Consider the sum
$$X(a_1\dots a_i,0) + \sum_{j=1}^{k-1}X(a_1\dots a_{\hati-1},j)$$
In the present case, $\hati\in\lgL(a)$. Our choice of $L\geq k$ is
now used, for it implies that $a_{\hati}\geq L\geq k$, and hence
$$X(a_1\dots a_i,0) + \sum_{j=1}^{k-1}X(a_1\dots a_{\hati-1},j)\leq X(a_1\dots a_i,0) + \sum_{j=1}^{a_{\hati}-1}X(a_1\dots a_{\hati-1},j)\periodeq$$
The terms in the above sum each appear once in
\begin{eqnarray} & & \sum_{t=i}^{m-1}\sum_{j=0}^{a_{t+1}-1}X(a_1\dots a_t,j) + X(a,0) +
\sum_{j=1}^{k-1}\sR_j(a) \\ \nonumber &=& \B(a)- \B(a_1\dots a_i)+
X(a,0) + \sum_{j=1}^{k-1}\sR_j(a)\end{eqnarray} and it follows
that
\begin{eqnarray}& & X(a_1\dots a_i,0) + \sum_{j=1}^{k-1}X(a_1\dots
a_{\hati-1},j)\\ \nonumber &\leq& \B(a)- \B(a_1\dots a_i)+ X(a,0)
+ \sum_{j=1}^{k-1}\sR_j(a)\periodeq\end{eqnarray}  Therefore,
\begin{multline}\label{eq:boundH11}F^T_a\text{  occurs  }\Rightarrow\B(a) + X(a,0) +\sum_{j=1}^{k-1}\sR_j(a)\leq \s \\ \Rightarrow X(a_1\dots a_i,0) + \sum_{j=1}^{k-1}X(a_1\dots a_{\hati-1},j)
\leq \s - \B(a_1\dots a_i)\periodeq\end{multline} But it is always
true that
$$\s \leq \B(a_1\dots a_{i-1}) +\p(a_1\dots a_{i-1})$$
and
$$\B(a_1\dots a_{i-1}) +\p(a_1\dots a_{i-1}) - \B(a_1\dots a_i) = \sum_{j\geq a_i}X(a_1\dots a_{i-1},j)\periodeq$$
Hence \begin{eqnarray} & & F^T_a\text{  occurs  }\\ \nonumber
&\Rightarrow & X(a_1\dots a_i,0) + \sum_{j=1}^{k-1}X(a_1\dots
a_{\hati-1},j) \leq \sum_{j\geq a_i}X(a_1\dots a_{i-1},j)\\
\nonumber &\Rightarrow& H_{a,i}\text{
occurs}\periodeq\end{eqnarray}

\mypar{Third containment: $H_{a,i}\subset F^T_a$ if $i\in\lgL(a)$
and $\hati=m+1$. } In this case, the terms of the sum
$$X(a_1\dots a_i,0) + \sum_{j=1}^{k-1}\sR_j(a)$$
are all contained in
\begin{eqnarray} & &\sum_{t=i}^{m-1}\sum_{j=0}^{a_{t+1}-1}X(a_1\dots a_t,j) + X(a,0) +
\sum_{j=1}^{k-1}\sR_j(a) \\ \nonumber&=& \B(a)-\B(a_1\dots a_i)
+X(a,0) + \sum_{j=1}^{k-1}\sR_j(a)\periodeq\end{eqnarray} The rest
of the proof proceeds exactly as in the case of the second
containment.\\

We now show that the events in \eqnref{KRmain1_Gagain} to
\eqnref{KRmain1_Hagain2} are {\em independent}. This is proven by
showing that no term $X(b,r)$ appears in the definition of more than one of those events. We will analyze three different cases.\\

\mypar{Comparing $G^T_a$ to the remaining events. }$G^T_a$ is
entirely defined in terms of $X(a_1\dots a_t,0)$ for
$t\in\smL(a)$. The only terms of the form $X(b,0)$ appearing in
the definition of the events $H_{a,i}$ have $b=a_1\dots a_i$ for
$i\in\lgL(a)$. This implies that no random variable appears in the
definition of both $G^T_a$ and $H_{a,i}$, for all $i\in\lgL(a)$.\\

\mypar{Comparing $H_{a,i}$ to $H_{a,\ell}$ for $i<\hati<\ell$,
$i,\ell\in  \lgL(a)$. } The definition of $H_{a,i}$ only involves
random variables of the form $X(a_1\dots a_t,j)$ for some
$j\in\N\cup\{0\}$ and $t\leq \hati-1<\ell-1$, whereas the
definition of $H_{a,\ell}$ involves $X(a_1\dots a_s,j)$ for $s\geq
\ell-1$. Therefore, the ranges of the indices $t$ and $s$ will
never overlap in this case.\\

 \mypar{Comparing $H_{a,i}$ to $H_{a,\ell}$ for $i<\hati=\ell$, $i,\ell\in\lgL(a)$.} By the same argument and with the same notation as above, the only "possibility for trouble is when $t=\hati-1=\ell-1=s$. This is precisely where the assumption that $L\geq k$ comes in. The event $H_{a,i}$ involves random variables of the form
 \begin{equation}\{X(a_1\dots a_{\hati-1},j)\suchthat 1\leq j\leq k-1\}\commaeq\end{equation}
 whereas the event $H_{a,\ell}$ uses the random variables
\begin{equation}\{X(a_1\dots a_{\ell-1},j)\suchthat j\geq a_\ell\}\periodeq\end{equation}
Since $\ell\in\lgL(a)$, $a_\ell>L\geq k$, the ranges of $j$ in the
two formulae above do not overlap, and we are done.\end{proof}

\subsection{Proof of \lemref{KRmain1_main1}}\label{sub:KRmain1_final}

Having proven \lemref{KRmain1_mixed}, we now come to the end of the proof of \lemref{KRmain1_main1}.\\

\begin{proof}{[of \lemref{KRmain1_main1}]} Our aim is to show that for any $T>0$,
\begin{equation}\label{eq:aimbnd}\sum_{a\in\Nstr}\Pr{a\mbox{ is $k$-fertile}\mid \p(\estr)\leq T}<+\infty\periodeq\end{equation}
To this end, we employ \lemref{KRmain1_mixed} and prove instead
that for some fixed number $L\geq L_0$ depending only on $f$, $k$
and $T$,
\begin{equation}\label{eq:aimbnd2}\sum_{a\in\Nstr} \frac{\max\{C,T\}^{m}}{|\smL(a)|!} \,\prod_{j\in\lgL(a)} \frac{1}{a_j^{(p-1)k}}<+\infty\periodeq\end{equation}
We will eventually choose some $L$ such that
\begin{equation}\label{eq:aimbndm}\sum_{|a|=m} \frac{\max\{C,T\}^{m}}{|\smL(a)|!} \,\prod_{j\in\lgL(a)} \frac{1}{a_j^{(p-1)k}} =2^{-\ohmega{m}} \mbox{ as }m\to +\infty\commaeq\end{equation}
which clearly implies \eqnref{aimbnd2}. Fix some $m$ and a subset
$S\subseteq \{1,\dots, m\}$ of size $|S|=s$. The sum of the above
quantities over all $a$ of length $|a|=m$ with $\smL(a)=S$  is
\begin{equation*}L^s
\frac{\max\{T,C\}^{m}}{s!}\prod_{i\in\{1,\dots,m\}\backslash
S}\sum_{a_i>L} \frac{1}{a_i^{(p-1)k}}\commaeq\end{equation*}
because there are $L^s$ ways of choosing the $a_j$'s with
$j\in\smL(a)$. Now note that $S$ can be chosen in $\binom{m}{s}$
for any $0\leq s\leq m$, and therefore
\begin{eqnarray}\label{eq:aimbndm2} & & \sum\limits_{\stackrel{|a|=m}{|\smL(a)|=s}}  \frac{{\max\{C,T\}}^{m}}{|\smL(a)|!} \,\prod_{j\in\lgL(a)} \frac{1}{a_j^{(p-1)k}}  \\ \nonumber &= &\binom{m}{s}L^s\frac{{\max\{C,T\}}^{m}}{s!} \,\prod_{j=1}^{m-s} \sum_{a_j> L}\frac{1}{a_j^{(p-1)k}} \\
\nonumber & \leq &\binom{m}{s}L^{s}\frac{{\hat T}^{m}}{s!}\left(\int_{L}^{+\infty}\frac{dx}{x^{k(p-1)}}\right)^{m-s} \\
\nonumber & \leq & \binom{m}{s}(L)^{s-\alpha(m-s)}\frac{{\hat
T}^{m}}{s!}\commaeq\end{eqnarray} where $\alpha\equiv k(p-1)-1$
and
$$\hat T\equiv \max\{C,T\}\times \max\left\{1,\frac{1}{\alpha}\right\}\periodeq$$
Here we make critical use of the condition $p>p_k=1+1/k$: under
this assumption, $\alpha>0$. Summing over $s$, we discover that
\begin{equation}\label{eq:aimbndsum}\sum_{|a|=m}\frac{\max\{C,T\}^{m}}{|\smL(a)|!} \,\prod_{j\in\lgL(a)} \frac{1}{a_j^{(p-1)k}} \leq \sum_{s=0}^m \binom{m}{s}\frac{\hat T^m L^{s-(m-s)\alpha}}{s!}\periodeq\end{equation}
To bound this last sum, we split it into two parts, corresponding
to $s\leq \alpha m/2(1+\alpha)$ and $s>\alpha m/2(1+\alpha)$. For
the first part, we forget the $s!$ term and bound
$s-(m-s)\alpha\leq -\alpha m/2$; for the second, we simply bound
$s!\geq \lceil\alpha m/2(1+\alpha)\rceil!$ and $s-(m-s)\alpha\leq
m$.
\begin{eqnarray}\nonumber \sum_{s\leq \frac{\alpha m}{2(1+\alpha)}} \binom{m}{s}\frac{\hat T^m L^{s-(m-s)\alpha}}{s!} &\leq &\left(\frac{\hat T}{L^{\frac{\alpha}{2}}}\right)^m \sum_{s\leq \frac{\alpha m}{2(1+\alpha)}} \binom{m}{s}\commaeq\\ \nonumber
\sum_{s> \frac{\alpha m}{2(1+\alpha)}} \binom{m}{s}\frac{\hat T^m
L^{s-(m-s)\alpha}}{s!} & \leq & \frac{\bigl(\hat
TL\bigr)^m}{\left\lceil\frac{\alpha m}{2(1+\alpha)}\right\rceil!}
\sum_{s> \frac{\alpha m}{2(1+\alpha)}}
\binom{m}{s}\periodeq\end{eqnarray} It follows that for $L \geq
(4\hat{T})^{2/\alpha}$, which only depends on $f$, $p$ and $T$,
\begin{eqnarray}\sum_{|a|=m} \binom{m}{s}\frac{\hat{T}^m L^{s-(m-s)\alpha}}{s!} &\leq &2^m \left( \frac{1}{4^m} + \frac{\bigl(\hat T^{1+2/\alpha}\bigr)^m}{\left\lceil\frac{\alpha m}{2(1+\alpha)}\right\rceil!} \right) \\ &=& 2^{-\ohmega{m}}\mbox{ as }m\to +\infty\periodeq\end{eqnarray}
This proves \eqnref{aimbndm} and finishes the proof.\end{proof}

\section{The structure of the infinite tree}\label{sec:KRmain2}
Now that the proof of \thmref{KRintro_main1} is complete, we
proceed to prove \thmref{KRintro_main2}. We will assume throughout
the section that $f(x)=(x+1)^p$ (with $p>1$) and that $k=k_p$ is
as in the statement of the Theorem. As in the previous section, it
is convenient to break the proof down into steps.

\begin{lemma}\label{lem:finchild}If $S$ is a rooted tree with $|S|=\ell+1$ vertices, then for all
$a\in\Nstr$
\begin{equation}\label{eq:finchild}\Pr{\#\{n\in\nat\; :\;
\W_{an}(\p(a)+\B(a)-\B(an))\mbox{ is isomorphic to }S\}
=+\infty}=\left\{\begin{array}{ll}1&p \leq p_\ell\commaeq\\ 0 &
p>p_\ell\periodeq\end{array}\right.\end{equation} \end{lemma}

\begin{lemma}\label{lem:posproblim}Let $\tilde \tree$ be any finite tree and $v$ be a vertex of $\tilde \tree$. There is a positive probability that all of the following events hold: \begin{enumerate}
\item the labelled $\KR$ process reaches state $\tilde\tree$;
\item $v$ is the unique vertex present in $\tilde\tree$ to have
any children after state $\tilde\tree$ is reached; and \item all
nodes that are born after state $\tilde\tree$ is reached are
$\ell$-fertile for some $\ell<k$.\end{enumerate}\end{lemma}

As we shall see below, these lemmas permit that
\thmref{KRintro_main2} is easily proven.\\

\begin{proof}{[of \thmref{KRintro_main2}]} By \claimref{explodesright}, there almost surely exists a unique node $v\in\Nstr$ with $\B(v)+\p(v)=\s$, and all other nodes have finitely many descendants in $\ftree$. Moreover, since $p>1+1/k_p$, one can apply \thmref{KRintro_main1} and deduce that with probability $1$ there are only finitely many children $vn$ of $v$ that are $k$-fertile.\\

If we remove all {\em other} children of $v$ (i.e. those that have
$\leq k-1$ descendants, which must be infinitely many) and their
descendants from $\ftree$, we obtain a finite tree $\tree$. We
{\em claim} that in fact $\ftree=\glue$. For consider some
(rooted, oriented) tree $S$ with $|S|\leq k$. by
\lemref{finchild}, there almost surely exist infinitely many
$n\in\nat$ such that $\W_{vn}(\p(v)-\sum_{j=0}^{n-1}X(v,j))$ is
isomorphic to $S$, and because $\s=\B(v)+\p(v) = \B(vn) +
(\p(v)-\sum_{j=0}^{n-1}X(v,j))$, this implies that
$\W_{vn}(\s-\B(vn))$ is isomorphic to $S$ for infinitely many $n$.
But $\W_{vn}(\s-\B(vn))$ is the subtree of $\ftree=\W(\s)$ rooted
at (and oriented towards) $vn$, hence with probability $1$ there
are infinitely many $n\in\nat$ such that the subtree of $\ftree$
rooted at $vn$ is isomorphic to $S$. This is true for any $S$ of
size $\leq k$, so all such trees must appear infinitely often, and finishes the proof of the claim.\\

We have shown that $\ftree$ is always isomorphic to some $\glue$.
Moreover, \lemref{posproblim} says that {\em any} $\glue$ has a
positive probability of being the value of $\ftree$. This finishes
the proof.\end{proof}

We now proceed to prove to prove Lemmas \ref{lem:finchild} and
\ref{lem:posproblim}.

\subsection{Proof of \lemref{finchild}}

\begin{proof}[of \lemref{finchild}] The $p>p_\ell$ case is implied by \thmref{KRintro_main1}, so we focus on $p\leq
p_\ell$, using ``$\approx$" to denote a rooted oriented tree
isomorphism. We will prove the theorem only for the case
$a=\estr$. This entails no loss of generality because the joint
distribution $\W_a(\cdot),\p(a)$ and $\{\B(an)-\B(a)\}_{n\in\nat}$
does
not depend on the choice of $a\in\Nstr$.\\

Define the sequence of events
\begin{equation}\label{eq:defAhere}B_n \equiv \{\W_n(\p(\estr)-\B(n))\approx S\}\;\;(n\in\nat)\periodeq\end{equation}
Our goal is to show that
$$\Pr{B_n\text{ infinitely often}}=1\periodeq$$
If the events $B_n$ were independent, we could apply the
Borel-Cantelli Lemma for independent events to prove this
statement. Since independence is lacking, we will substitute the
events $B_n$ by a sequence of independent events $A_n$ such that
\begin{eqnarray}\label{eq:cansubs}\Pr{A_n\text{ i.o. but not }B_n\text{ i.o.}} & = & 0\periodeq\\
\label{eq:infsum} \sum_{n\in\nat}\Pr{A_n}
&=&+\infty\periodeq\end{eqnarray} Because the sequence $\{A_n\}$
consists of independent events, equation \eqnref{infsum} implies
that $A_n$ infinitely often almost surely, which implies (via
equation \eqnref{cansubs}) that $B_n$ infinitely often almost
surely. Therefore, \eqnref{cansubs} and \eqnref{infsum} imply the
Lemma.\\

We define the sequence $A_n$ as follows
\begin{equation}A_n \equiv \left\{\forall t\in\left[\frac{1}{2(p-1)n^{p-1}},\frac{3}{(p-1)n^{p-1}}\right] \;\;\W_n(t)\approx S\right\}\;\;(n\in\nat)\periodeq\end{equation}
The independence of those events is a consequence of the
independence of the processes $\{\W_n(\cdot)\}_{n\in\nat}$.
Moreover,
\begin{multline}\label{eq:BCcool}\Pr{A_n\text{ i.o. but not }B_n\text{ i.o.}}\\ \leq
\Pr{\p(\estr)-\B(n)\not\in\left[\frac{1}{2(p-1)n^{p-1}},\frac{3}{2(p-1)n^{p-1}}\right]\text{
i.o.}}\periodeq\end{multline} We {\em claim} that the event on the
RHS of \eqnref{BCcool} has probability $0$. To see this, note that
$$\p(\estr)-\B(n) = \sum_{j\geq n}X(\estr,j)$$
is a sum of independent, rate-$f(j)$ exponentials, and
$$\Ex{\p(\estr)-\B(n)} = \sum_{j\geq n}\frac{1}{f(j)}\periodeq$$
As a result, direct use of \lemref{lgdev} and the estimate
$$S_1(n) = \sum_{j=n}^{+\infty}\frac{1}{(j+1)^p} \sim \frac{1}{(p-1)n^{p-1}}\;\;\;(n\gg 1)\commaeq$$
implies
$$\sum_{j\geq 1}\Pr{\p(\estr)-\B(n)\not\in\left[\frac{1}{2(p-1)n^{p-1}},\frac{3}{2(p-1)n^{p-1}}\right]} <+\infty\periodeq$$
Therefore, the Borel-Cantelli Lemma implies that
\begin{equation}\label{eq:bcconv}\Pr{\p(\estr)-\B(n)\not\in\left[\frac{1}{2(p-1)n^{p-1}},\frac{3}{2(p-1)n^{p-1}}\right]\text{
i.o.}}=0\commaeq\end{equation} thereby proving the claim and (via
\eqnref{BCcool}) equation \eqnref{cansubs}.\\

It remains to prove \eqnref{infsum}. For this purpose, we will
only need a very rough lower bound on the probability of $A_n$.
Consider a {\em labelling} of the elements of $S$. That is, pick a
finite \proper~ subset of $\Nstr$, i.e. an element $\hat S\in
\Efin$, that corresponds to a labelling of the vertex set of $S$
as defined in \subref{KRintro_label}. We assume that $\hat S$ is
ordered
\begin{equation}\label{eq:labord}\hat S = \{s^{(0)}=\estr,s^{(1)},\dots s^{(\ell)}\}\end{equation}
in a way such that for all $1\leq i\leq \ell$, there is an index
$p_i<i$ such that $s^{(p_i)}$ is the parent sequence of $s^{(i)}$.
We also define the subsets
$$\hat S(i)\equiv \{s^{(0)},\dots,s^{(i)}\}\;\;(1\leq i\leq \ell)$$
The ordering property implies that $\hat S(i)$ is also a \proper~
subset of $\Nstr$. Now define (for $1\leq i\leq \ell$, where
applicable):
\begin{eqnarray}\label{eq:tn} t_n &=& \frac{1}{2(p-1)n^{p-1}}\commaeq \\ \label{eq:Tn}
T_n & = & \frac{3}{2(p-1)n^{p-1}}\commaeq \\ \label{eq:defCni}
C_n(i)& = & \left\{\text{
$s^{(i)}$ is the only vertex born in }\{\W(t)\}_{t\in \left[\frac{i-1}{\ell}t_n,\frac{i}{\ell}t_n\right]}\right\}\commaeq\\
\label{eq:defDn} D_n &=& \left\{\text{no vertex is born in
}\{\W(t)\}_{t\in[t_n,T_n]}\right\}\periodeq\end{eqnarray} Clearly,
\begin{equation}\label{eq:lbcap}\Pr{A_n} \geq \Pr{C_n(1)}\times \left\{\prod_{i=2}^{\ell}\Pr{C_n(i)\left| \bigcap_{j=1}^{i-1}C_n(j)\right.} \right\}\times \Pr{D_n\left|\bigcap_{r=1}^{\ell}C_n(r)\right.}\periodeq\end{equation}
(In fact, $A_n$ is defined in terms of $\W_n(\cdot)$ rather than
$\W(\cdot)$, but in terms of evaluating the probabilities that
does not make any difference since these two processes have the
same distribution.) We will lower bound the probabilities on the
RHS of the above inequality. \\

\mypar{Probability of $C_n(1)$. }The probability of $C_n(1)$ is
the probability that the birth time of $s^{(1)}$ is
$\B(s^{(1)})\leq t_n/\ell$ and that no other birth occurs in the
time interval $[\B(s^{(1)}),t_n/\ell]$. Conditioning on a value
$0\leq \B(s^{(1)})= t\leq t_n/\ell$, the time of the next birth in
$\W(\cdot)$ is
$$\min\{X(s^{(0)},0),X(s^{(1)},1)\} \eqdist \exp(f(0)+f(1))\periodeq$$
Hence \begin{multline}\Pr{C_n(1)\mid \B(s^{(1)})=u } =
\Pr{\exp(f(0)+f(1)) \geq \frac{t_n}{\ell}-u} \\ \geq
e^{-(f(0)+f(1))\frac{t_n}{\ell}-u}\geq
e^{-(2f(2))(\frac{t_n}{\ell})} \periodeq\end{multline} Moreover,
$$\Pr{\B(s^{(0)})\leq \frac{t_n}{\ell}} = 1 - e^{-f(0)\frac{t_n}{\ell}}\periodeq$$
Since $t_n\to 0$ as $n\to +\infty$, it follows that there exist
constants $C_1,n_1>0$ such that for all $n\geq n_1$
$$\Pr{\B(s^{(0)})\leq \frac{t_n}{\ell}}\geq C_1 t_n\periodeq$$
We conclude that
\begin{equation}\label{eq:prc1}\Pr{C_n(1)}\geq C_1 t_n\, e^{-2f(2)\frac{t_n}{\ell}}\periodeq\end{equation}

\mypar{Probability of $C_n(i)$, $2\leq i\leq \ell$. }In this part,
we will make use of the Markov property of the continuous-time
process (cf. \remref{CTMC}). Notice that the conditioned event is
defined entirely in terms of $\{\W(s)\}_{0\leq s\leq
(i-1)t_n/\ell}$, whereas $C_n(i)$ is defined entirely in terms of
$\{\W(s)\}_{s\geq (i-1)t_n/\ell}$. Moreover, it is also true that
inside the event $\cap_{j=1}^{i-1}C_n(j)$
\begin{equation}\label{eq:cond}\W\left(\frac{i-1}{\ell}t_n\right) = \hat
S(i-1)\periodeq\end{equation} Therefore, we can apply
\remref{CTMC} to deduce
\begin{equation}\Pr{C_n(i)\left| \bigcap_{j=1}^{i-1}C_n(j)\right.} = \Pr{C_n(i)\left| \W\left(\frac{i-1}{\ell}t_n\right) = \hat
S(i-1)\right.}\periodeq\end{equation} For $C_n(i)$ to happen, two
conditions must be satisfied.
\begin{enumerate}
\item $\B(s^{(i)})-({i-1})t_n/\ell \leq {t_n}{\ell}$. That is,
$s^{(i)}$ must be born in the interval
$$[(i-1)t_n/\ell,{i}t_n/\ell]\periodeq$$ \item No other birth
happens in the interval $[({i-1})t_n/\ell,it_n/\ell]$.
\end{enumerate}
Choose a value $0\leq u\leq t_n/\ell$. We will now bound
\begin{equation}\label{eq:condmore}\Pr{C_n(i)\left|\W\left(\frac{(i-1)t_n}{\ell}\right) = \hat
S(i-1),\;\B(s^{(i)})=u+\frac{(i-1)t_n}{\ell}\right.}\periodeq\end{equation}
In this case, note that the rate at which the first birth of a
node $a\neq s^{(i)}$ happens in $\hat{S}(i-1)$ is\footnote{This is
the rate {\em until some birth happens}, whether it is the birth
of $a$ or of some $s^{(i)}\neq a$.}
\begin{equation}R_{i-1}\equiv \sum_{0\leq j<i,\; j\neq p_i}f(\deg_{\hat{S}(i-1)}(s^{(j)})) \leq (i+1)f(i+1)\commaeq\end{equation}
the inequality being justified  by the fact that the cardinality
of $\hat{S}(i-1)$ is $i$. The rate of births {\em after} time
$(i-1)t_n/\ell+u$ under the conditioning of \eqnref{condmore} is
\begin{equation}T_{i}\equiv \sum_{0\leq j\leq i}f(\deg_{\hat{S}(i)}(s^{(j)})) \leq (i+1)f(i+1)\periodeq\end{equation}
Under the conditioning in \eqnref{condmore}, $C_n(i)$ holds iff no
$a\neq s^{(i)}$ is born in the time interval
$[(i-1)t_n/\ell,(i-1)t_n/\ell+u]$ and no births happen in
$[(i-1)t_n/\ell+u,it_n/\ell]$. By the Markov property of
$\W(\cdot)$, these events in different time intervals are
independent given $\W((i-1)t_n/\ell+u)=\hat{S}(i)$. Therefore, we
can write
\begin{multline}\Pr{C_n(i)\left|\W\left(\frac{(i-1)t_n}{\ell}\right)
= \hat S(i-1),\;\B(s^{(i)})=u+\frac{(i-1)t_n}{\ell}\right.} \\ =
\Pr{\exp(R_{i-1})>u} \times \Pr{\exp(T_i)>t_n/\ell-u} \\ =
e^{-{R_{i-1}u}}e^{-{T_{i}}(\frac{t_n}{\ell}-u)}\geq
e^{-(i+1)f(i+1)\frac{t_n}{\ell}}\periodeq\end{multline} As a
result,
\begin{multline}\Pr{C_n(i)\left|\W\left(\frac{(i-1)t_n}{\ell}\right)
= \hat S(i-1)\right.} \\ \geq e^{-(i+1)f(i+1)\frac{t_n}{\ell}}
\Pr{\B(s^{(i)})-\frac{(i-1)t_n}{\ell}\leq
\frac{t_n}{\ell}\left|\W\left(\frac{(i-1)t_n}{\ell}\right) = \hat
S(i-1)\right.}\periodeq\end{multline} Now notice that conditioned
on $\W((i-1)t_n/\ell)$,
$$\B(s^{(i)})-\frac{(i-1)t_n}{\ell}\eqdist\exp(f(\deg_{\hat{S}(i-1)}(s^{(p_i)})))\periodeq$$
Hence
\begin{multline}\Pr{\B(s^{(i)})-\frac{(i-1)t_n}{\ell}\leq
\frac{t_n}{\ell}\left|\W\left(\frac{(i-1)t_n}{\ell}\right)= \hat
S(i-1)\right.} \\ = 1 -
e^{-f(\deg_{\hat{S}(i-1)}(s^{(p_i)}))\frac{t_n}{\ell}}\geq 1 -
e^{-f(\ell)\frac{t_n}{\ell}}\periodeq\end{multline} To state our
bound for the probability $C_n(i)$, we note that $t_n\to 0$ as
$n\to +\infty$, and therefore there exist constants $C_i,n_i>0$
depending only on $\ell$ and $f$ such that for all $n\geq n_i$
\begin{equation}\label{eq:finallyCi}\Pr{C_n(i)\left|\W\left(\frac{(i-1)t_n}{\ell}\right)
= \hat S(i-1)\right.}\geq e^{-(i+1)f(i+1)\frac{t_n}{\ell}} C_i
t_n\periodeq\end{equation}

\mypar{Probability of $D_n$. }For this bound, we again use the
Markov property of $\W(\cdot)$. Notice that whereas $D_n$ is only
defined in terms of $\{\W(t)\}_{t\geq t_n}$, the definition of
$\cup_{i\leq \ell}C_n(i)$ only depends on $\{\W(t)\}_{0\leq t\leq
t_n}$. Moreover, inside the latter event, $\W(t_n) = \hat S$. We
can then apply \remref{CTMC} to conclude
\begin{equation}\Pr{D_n\left|\bigcup_{i=1}^\ell C_n(i)\right.} = \Pr{D_n\left|\W(t_n)=\hat S\right.}\periodeq\end{equation}
Under this last conditioning, the rate of new births in $\W(t_n)$
is
$$\sum_{s\in S}f(\deg_S(s)) \leq (\ell+1)f(\ell+1)$$
and the probability that none of those births occur in $[t_n,T_n]$
is precisely
\begin{equation}\label{eq:finallyD}\Pr{D_n\left|\W(t_n)=\hat S\right.} = e^{-(T_n-t_n)\sum_{s\in S}f(\deg_S(s)) }\geq e^{-(T_n-t_n)(\ell+1)f(\ell+1) }\periodeq\end{equation}
\mypar{Wrapping up. }To finish this proof, we plug \eqnref{prc1},
\eqnref{finallyCi} and \eqnref{finallyD} into \eqnref{lbcap},
letting $n\geq \max\{n_i:1\leq i\leq \ell\}$ and $C=C_1C_2\dots
C_\ell$
\begin{eqnarray}\label{eq:therefinprob}\Pr{A_n}&\geq& C\left(\prod_{i=1}^\ell e^{-(i+1)f(i+1)\frac{t_n}{\ell}}\right)\times e^{-(T_n-t_n)(\ell+1)f(\ell+1)} \times t_n^{\ell}\\
&\geq & C \exp(-(\ell+1)f(\ell+1)T_n)
t_n^{\ell}\periodeq\end{eqnarray} For $\ell$ fixed, $n\to
+\infty$, we deduce (using the definition of $t_n$ and $T_n$ in
\eqnref{tn}, \eqnref{Tn})
$$\Pr{A_n}= \ohmega{t_n^{\ell}} = \ohmega{n^{-(p-1)\ell}}\periodeq$$ Now the assumption $p\leq p_\ell$ comes
into play, for it implies that $(p-1)\ell\leq 1$. As a result
$$\sum_{n\geq 1}\Pr{A_n} = \sum_{n\geq 1} \ohmega{n^{-(p-1)\ell}} = +\infty\periodeq$$
This proves \eqnref{infsum} and finishes the proof.\end{proof}

\subsection{Proof of \lemref{posproblim}}

\begin{proof}{[of \lemref{posproblim}]} Let $\tree$ have $u+1$ vertices and let $v$ have $r$ children in $\tree$.  Let $(\tree,v)+n$ denote $\tree$ with $n$ additional children added
to $v$.  Asymptotically in $n$ we consider the probability that
$\ftree$ is not isomorphic to $\glue$ conditional on the $\KR$
process reaching $(\tree,v)+n$.  Each $w\in \tree$, $w\neq v$, has
probability $o(1)$ of having a child before $v$ explodes.  Each of
the $n$ additional children of $v$ has probability $O(n^{k(1-p)})$
of having $k$ (or more) descendants before $v$ explodes.  For
$i>n+l$ the $i$-th child of $v$ has $k$ (or more) descendants
before $v$ explodes with probability $O(i^{k(1-p)})$.  The total
probability of any of these events occurs is then bounded from
above by $u\cdot o(1) + n\cdot O(n^{k(1-p)}) +
\sum_{i>n+l}O(i^{k(1-p)})$ which is $o(1)$ because $k(1-p)<-1$. We
can therefore find an explicit $n$ so that this probability is
less than, say,
$\frac{1}{2}$.\\

With positive (perhaps small) probability the first $n+l$ steps of
the $\KR$ process yield $(T,v)+n$.  Then with probability at least
$\frac{1}{2}$ the final $T_{\infty}$ is $\glue$ as
desired.\end{proof}

\section{Conclusion}\label{sec:KRextra}
The two main theorems of this paper completely characterize the
limits of the super-linear $\KR$ process. Some of their
consequences are the fact that the tree $\ftree$ has finite height
(and thus the finite-time $\KR$ trees have bounded height), and
that the nodes of in-degree $\leq k-1$ (where $p>p_k$) are all but
finitely many. However, these characteristics raise many
interesting questions about {\em distributions} of the above
quantities. For instance, what does the tail of the height
distribution of $\ftree$ look like? We believe that the methods
presented in this paper might be sharpened to prove this and other
results.\\

There are many more open questions about the $p<1$ case of $\KR$.
The authors of \cite{Chung03} have derived some results on their
modified model for this range of $p$ under the assumption that
certain limits exist. Proving unconditional results of this nature
for the $\KR$ model remains an important open problem that is also
potentially amendable to treatment by our techniques, since the
exponential embedding applies to any attachment kernel.\\

It would also be quite interesting if the exponential embedding
could be used to prove known and new properties of related network
models, in particular the original Barab\' asi-Albert preferential
attachment model. The rigorous version of the process defined in
\cite{Bollobas04} is essentially the $\KR$ process defined in our
paper with attachment kernel $f(x)=x+1$, and it could be the case
that the embedding method is a viable technical alternative to the
``linearized chord diagrams" of \cite{Bollobas04}.

\appendix

\section{Appendix -- proofs of technical lemmas}\label{appendix}
\begin{proof}[of \lemref{lgdev}]  We will only prove the first
inequality, for the proof of the second one is very similar.  The
technique we employ is fairly standard and is commonly used in
other proofs of Chernoff-type large deviation inequalities
\cite{AlonSpencer_Method}. Let $A_n = \sum_{j\geq n}X_j -
f(j)^{-1}$. Fix any $0<s\leq (n+1)^p/2$ and notice that, by the
standard Bernstein's trick, the formulae in \subref{probing}, the
inequality ``$1+x\leq e^x$'', and some simple calculations
\begin{eqnarray*}\Pr{A_n\,>\,\delta} &=& \Pr{e^{s\,A_n}>e^{s\,\delta}}\\
& \leq &  e^{-s\,\delta}\Ex{e^{\sum_{j\geq n}s \left(X_j - \frac{1}{f(j)}\right)}} \\
& = &  e^{-s\,\delta}\prod_{j\geq n}\Ex{e^{s \left(X_j - \frac{1}{f(j)}\right)}}\\
& = &  e^{-s\,\delta}\prod_{j\geq n}\frac{e^{-\frac{s}{f(j)}}}{1-\frac{s}{f(j)}} \\
& = &  e^{-s\,\delta} \times \\ && \prod_{j\geq
n}e^{-\frac{s}{f(j)}}\left(1+
\frac{s}{f(j)} + \frac{s^2}{f(j)^2}\frac{1}{1-\frac{s}{f(j)}}\right) \\
&\leq & e^{-s\,\delta}\prod_{j\geq n}\exp(2\frac{s^2}{(j+1)^{2p}})\\
\label{eq:est_lastexp}&\leq& \exp\left(\frac{2s^2}{(2p-1)n^{2p-1}}
- s\delta\right)\end{eqnarray*} To finish the proof, we set
$s\equiv n^{p-1/2}$, which is permissible since $n^{p-1/2}\leq
(n+1)^p/2$ for all large enough $n$.\end{proof}

\begin{proof}{[of Lemma $\ref{lem:prprod}$]}To begin with, we note that
\begin{equation*}\label{eq:WWWWWWW}\forall \eps>0 \; \prod_{i=1}^k \pr{Y_i\leq \frac{\eps}{k}}
\leq \pr{Y\leq \eps} \leq \prod_{i=1}^k \pr{Y_i\leq
\eps}\end{equation*} and therefore the assumptions imply the
existence of a constant $C_0$ depending only on the distributions
of the $Y_i$'s and on $k$ such that
\begin{equation*}\label{eq:aaaaauuuuu}\forall \eps>0\commaeq \;\;\;
\frac{\eps^k}{C_0} \leq \Pr{Y\leq \eps} \leq
C_0\eps^k\end{equation*} We now use the notation and results in
the proof of \lemref{lgdev} with $\delta = {n^{3/4-p}}$. Then
$$\forall n\geq n_0\commaeq\;\; \Pr{\left|\sum_{j\geq n}X_j - \mu \right|>\delta}\leq 2C_1\, e^{-n^{\frac{1}{4}}}
$$
for some constant $C_1$. Then
\begin{multline*} \Pr{Y\leq \mu-\delta} - \Pr{\left|\sum_{j\geq n}X_j - \mu \right|>\delta} \\ \leq \Pr{Y\leq \sum_{j\geq n}X_j} \leq
\Pr{Y\leq \mu+\delta} + \Pr{\left|\sum_{j\geq n}X_j - \mu
\right|>\delta}\commaeq \end{multline*} and by the previous bounds
$$ \frac{(\mu-\delta)^k}{C_0} - 2C_1\, e^{-n^{\frac{1}{4}}} \leq \Pr{Y\leq \sum_{j\geq n}X_j} \leq C_0(\mu+\delta)^k + 2C_1\, e^{-n^{\frac{1}{4}}}\periodeq$$
The result now follows from the fact that, as $n\to +\infty$
$$\mu \sim \frac{1}{(p-1)n^{p-1}} \sim \mu\pm\delta \gg e^{-n^{\frac{1}{4}}}\periodeq$$

\ignore{as $n\to +\infty$ \begin{eqnarray*}\mu\pm \delta &\sim&
\frac{1}{(p-1)n^{p-1}} \\ \Pr{\left|\sum_{j\geq n}X_j - \mu
\right|>\delta}&\ll& \frac{1}{n^{(p-1)k}}\end{eqnarray*} Equation
\eqnref{WWWWWWW} and the large-deviation bound imply
$$\Pr{Y\leq \mu\pm \delta} \pm \Pr{\left|\sum_{j\geq n}X_j - \mu
\right|<\delta}\sim \left(\frac{1}{(p-1)n^{p-1}}\right)^k$$ which
implies the lemma.}\end{proof}

\begin{proof}{[of Lemma $\ref{lem:Poi}$]} $\Pr{Z\leq \lambda}$ is equal to the probability that there are at least $k$
arrivals up to time $\lambda$ in a Poisson process with rate $1$.
This has a Poisson distribution with rate $\lambda$; hence we have
the exact result $$\Pr{Z\leq \lambda} =
\sum_{j=k}^{+\infty}e^{-\lambda}\frac{\lambda^{j}}{j!}$$ The upper
bound follows from $$ \sum_{j\geq k} \frac{\lambda^j}{j!}=
\frac{\lambda^k}{k!}\sum_{\ell\geq 0}
\frac{\lambda^{\ell}}{(k+1)(k+2)\dots (k+\ell)} \leq
\frac{\lambda^k e^\lambda}{k!}$$\end{proof}

\bibliography{connect_IM}   

\begin{thebibliography}{10}

\bibitem{AlbertSurvey}
R\'eka Albert and Albert-L\'aszl\'o Barab\'asi.
\newblock Statistical mechanics of complex networks.
\newblock {\em Reviews of Modern Physics}, 74:47--97, 2002.
\newblock Available at \texttt{cond-mat/0106096}.

\bibitem{Albert99}
R\'eka Albert, Haowoong Jeong, and Albert-L\'aszl\'o Barab\'asi.
\newblock Diameter of the {World Wide Web}.
\newblock {\em Nature}, 401:130--131, 1999.

\bibitem{AlonSpencer_Method}
Noga Alon and Joel Spencer.
\newblock {\em The Probabilistic Method}.
\newblock Wiley-Interscience Series in Discrete Mathematics. John Wiley and
  Sons, New York, second edition, 2000.

\bibitem{Barabasi99}
Albert-L\'aszl\'o Barab\'asi and R\'eka Albert.
\newblock Emergence of scaling in random networks.
\newblock {\em Science}, 286:509--512, 1999.

\bibitem{Bianconi01}
Ginestra Bianconi and Albert-L\'aszl\'o Barab\'asi.
\newblock Bose-{Einstein} condensation in complex networks.
\newblock {\em Physical Review Letters}, 86:5632–5635, 2001.

\bibitem{Bollobas03_2}
B\'ela Bollob\'as, Christian Borgs, Jennifer Chayes, and Oliver Riordan.
\newblock Directed scale-free graphs.
\newblock In {\em Proceedings of the 12th Annual {ACM-SIAM} Symposium on
  Discrete Algorithms}, pages 132--139. Society for Industrial and Applied
  Mathematics, Philadelphia, PA, USA, 2003.

\bibitem{BollobasSurvey}
B\'ela Bollob\'as and Oliver Riordan.
\newblock Mathematical results on scale-free random graphs.
\newblock In {\em Handbook of graphs and networks}, pages 1--34. Wiley-VCH,
  Weinheim.

\bibitem{Bollobas03}
B\'ela Bollob\'as and Oliver Riordan.
\newblock Robustness and vulnerability of scale-free random graphs.
\newblock {\em Internet Mathematics}, 1(1):1--35, 2003.

\bibitem{Bollobas04}
B\'ela Bollob\'as and Oliver Riordan.
\newblock The diameter of a scale-free random graph.
\newblock {\em Combinatorica}, 4:5--34, 2004.

\bibitem{Bollobas00}
B\'ela Bollob\'as, Oliver Riordan, Joel Spencer, and G\'abor T\'ardos.
\newblock The degree sequence of a scale-free random graph process.
\newblock {\em Random Structures and Algorithms}, 18(3):279--290, 2000.

\bibitem{Chung03}
Fan Chung, Shirin Handjani, and Doug Jungreis.
\newblock Generalizations of {P\'olya}'s urn problem.
\newblock {\em Annals of Combinatorics}, 7(2):141--153, 2003.

\bibitem{Cooper03}
Colin Cooper and Alan Frieze.
\newblock On a general model of web graphs.
\newblock {\em Random Structures and Algorithms}, 22:311--335.

\bibitem{Davis90}
Burgess Davis.
\newblock Reinforced random walk.
\newblock {\em Probability Theory and Related Fields}, 84(2):203--229, 1990.

\bibitem{MendesSurvey}
Serguei~N. Dorogovtsev and Jos\'{e}~F.F. Mendes.
\newblock Evolution of networks.
\newblock {\em Advances in Physics}, 51:1079--1187, 2002.
\newblock Available at \texttt{cond-mat/0106144}.

\bibitem{DrineaEM}
Eleni Drinea, Mihaela Enachescu, and Michael Mitzenmacher.
\newblock Variations on {Random Graph} models of the {Web}.
\newblock {Harvard Technical Report TR}-06-01, 2001.

\bibitem{Drinea02}
Eleni Drinea, Alan Frieze, and Michael Mitzenmacher.
\newblock Balls in bins processes with feedback.
\newblock In {\em Proceedings of the 11th Annual {ACM-SIAM} Symposium on
  Discrete Algorithms}, pages 308--315. Society for Industrial and Applied
  Mathematics, Philadelphia, PA, USA, 2002.

\bibitem{Khanin01}
Kostya Khanin and Raya Khanin.
\newblock A probabilistic model for the establishment of neuron polarity.
\newblock {\em Journal of Mathematical Biology}, 42(1):26--40, 2001.

\bibitem{Krapivsky02}
P.L. Krapivsky and Sidney~L. Redner.
\newblock Organization of growing random networks.
\newblock {\em Physics Reviews E}, 63:066123, 2001.
\newblock Available at \texttt{cond-mat/0011094}.

\bibitem{MitzenmacherOS04}
Michael Mitzenmacher, Roberto Oliveira, and Joel Spencer.
\newblock A scaling result for explosive processes.
\newblock {\em Electronic Journal of Combinatorics}, 11(1):R31, 2004.

\bibitem{Spencer??}
Joel Spencer and Nicholas Wormald.
\newblock Explosive processes.
\newblock Manuscript.

\end{thebibliography}

\end{document}